\documentclass[MSNbibl,number,citesort,seceqn,dvips]{arxbj}
\usepackage{mathbh,upgreek}

\aid{0}
\volume{19}
\issue{2}
\pubyear{2013}
\firstpage{646}
\lastpage{675}
\doi{10.3150/11-BEJ408}

\makeatletter
\newcommand{\R}{\mathbb{R}}
\newcommand{\N}{\mathbb{N}}
\newcommand{\E}{\mathbb{E}}
\newcommand{\Pro}{\mathbb{P}}
\newcommand{\eps}{\varepsilon}
\newcommand{\1}{\mathbh{1}}
\newcommand\cD{\mathcal{ D}}
\newcommand\cL{\mathcal{ L}}
\newcommand\cO{\mathcal{ O}}
\newcommand\cU{\mathcal{ U}}
\newcommand\cX{\mathcal{ X}}
\newcommand\cZ{\mathcal{ Z}}
\newcommand{\Z}{\mathbb{Z}}
\newcommand{\norm}[1]{\|#1\|}%
\renewcommand{\hat}{\widehat}
\newtheorem{Theorem}{Theorem}[section]
\newtheorem{theorema}{Theorem}
\newtheorem{Lemma}[Theorem]{Lemma}
\newproclaim{Definition}[Theorem]{Definition}
\newtheorem{Proposition}[Theorem]{Proposition}
\newtheorem{Corollary}[Theorem]{Corollary}
\newproclaim{Question}[Theorem]{Question}
\newcommand{\eqref}[1]{(\ref{#1})}
\makeatother

\begin{document}
\begin{frontmatter}

\title{On the optimality of the aggregate with exponential weights for low temperatures}
\runtitle{On the optimality of the aggregate with exponential weights}

\begin{aug}
\author[1]{\fnms{Guillaume} \snm{Lecu\'e}\corref{}\thanksref{1}\ead[label=e1]{guillaume.lecue@univ-mlv.fr}}%
\and
\author[2]{\fnms{Shahar} \snm{Mendelson}\thanksref{2}\ead[label=e2]{shahar@tx.technion.ac.il}}
\runauthor{G. Lecu\'e and S. Mendelson} 
\address[1]{CNRS, LAMA, Universit{\'e} Paris-Est Marne-la-vall{\'e}e,
Champs-sur-Marne
77454 France.\\ \printead{e1}}
\address[2]{Department of Mathematics, Technion, I.I.T, Haifa
32000, Israel.\\ \printead{e2}}
\end{aug}

\received{\smonth{11} \syear{2010}}
\revised{\smonth{7} \syear{2011}}

%
\begin{abstract}
Given a finite class of functions $F$, the problem of aggregation
is to construct a procedure with a risk as close as possible to
the risk of the best element in the class. A classical procedure
(PAC-Bayesian statistical learning theory (2004) Paris 6, \textit
{Statistical Learning
Theory and Stochastic Optimization} (2001) Springer, \textit{Ann.
Statist.} \textbf{28} (2000) 75--87) is the aggregate with
exponential weights (AEW), defined by
\[
\tilde{f}^{\mathrm{AEW}}=\sum_{f\in F}\hat\theta(f) f, \qquad
\mbox
{where } \hat\theta(f)=\frac{\exp(-({n}/{T})R_n(f))}
{\sum_{g\in F}\exp(-({n}/{T})R_n(g))},
\]
where $T>0$ is called the temperature parameter and $R_n(\cdot)$ is an
empirical risk.

In this article, we study the optimality of the AEW in the
regression model with random design and in the low-temperature
regime. We prove three properties of AEW. First, we show that AEW is a
suboptimal aggregation procedure in expectation with respect to
the quadratic risk when $T\leq c_1$, where $c_1$ is an absolute
positive constant (the low-temperature regime), and that it is
suboptimal in probability even for high temperatures. Second, we
show that as the cardinality of the dictionary grows, the behavior
of AEW might deteriorate, namely, that in the low-temperature
regime it might concentrate with high probability around elements
in the dictionary with risk greater than the risk of the best
function in the dictionary by at least an order of $1/\sqrt{n}$.
Third, we prove that if a geometric
condition on the dictionary (the so-called ``Bernstein condition'') is assumed,
then AEW is indeed optimal both in high probability and in
expectation in the low-temperature regime. Moreover, under that
assumption, the complexity term is essentially the logarithm of the
cardinality of the set of ``almost minimizers'' rather than the
logarithm of the cardinality of the entire dictionary. This result
holds for small values of the temperature parameter, thus
complementing an analogous result for high temperatures.
\end{abstract}

%

\end{frontmatter}

\section{Introduction and main results} \label{sec:introduction}
In this note we study the problem concerning the optimality of the AEW
in the regression model with random design. To formulate the problem,
we need to introduce several definitions.

Let $\cZ$ and $\cX$ be two measure spaces, and set $Z$ and
$Z_1,\ldots,Z_n$ to be $n+1$ i.i.d. random variables with values in
$\cZ$. From a statistical standpoint,
$\cD=(Z_1,\ldots,Z_n)$ is the set of given data at our disposal. The
\textit{risk}
of a measurable real-valued function $f$ defined on $\cX$ is given by
\[
R(f)=\E Q(Z,f),
\]
where $Q\dvtx\cZ\times\cL(\cX)\mapsto\R$ is a non-negative
function, called
the \textit{loss function} and $\cL(\cX)$ is the set of all
real-valued measurable functions defined on $\cX$. If $\hat f$ is a statistic
constructed using the data $\cD$, then the risk of $\hat f$ is the random
variable
\[
R(\hat f)=\E[Q(Z,\hat f)\vert\cD].
\]
Throughout this article, we restrict our attention to functions
$f$, loss functions $Q$, and random variables $Z$ for which
$|Q(Z,f)| \leq b$ almost surely. (Note that some results have been
obtained in the same setup for unbounded loss functions in
\cite{MR2163920,MR1946426,JRT08}, and~\cite{AoS08}.)
The loss function on which we focus throughout most of the article
is the quadratic loss function, defined when $Z=(X,Y)$ by
$Q((X,Y),f)=(Y-f(X))^2$.

In the aggregation framework, one is given a finite set $F$ of
real-valued functions defined on~$\cX$, usually called a
\textit{dictionary}. The problem of \textit{aggregation} (see, e.g.,
\cite{MR1775638,MR2163920}, and~\cite{MR1762904})
is to construct a procedure, usually called an \textit{aggregation
procedure}, that produces a function with a risk as close as
possible to the risk of the best element in $F$. Keeping this in
mind, one can define the \textit{optimal rate of aggregation}
\cite{TsyCOLT07,LM1}, which is the smallest price, as a function
of the cardinality of the dictionary $M$ and the sample size $n$,
that one has to pay to construct a function with a risk as close
as possible to that of the best element in the dictionary. We
recall the definition for the ``expectation case;'' a similar
definition for the ``probability case'' can be formulated
as well (see, e.g.,~\cite{LM1}).

\begin{Definition}[(\cite{TsyCOLT07})]\label{def:definition-optimality}
Let $b>0$. We say that $(\psi_n(M))_{n,M\in\N^*}$ is an optimal
rate of aggregation in expectation when there exist two positive
constants, $c_0$ and $c_1$, depending only on $b$, for which the
following holds for any $n\in\N^*$ and $M\in\N^*$:
\begin{enumerate}
\item There exists an aggregation procedure $\tilde f_n$ such that for
any dictionary $F$ of cardinality $M$ and any random variable $Z$
satisfying $|Q(Z,f)| \leq b$ almost surely for all $f\in F$, one
has
%
%
\begin{equation}\label{eq:exact-oracle-inequality}
\E R(\tilde f_n)\leq\min_{f\in F}R(f)+c_0\psi_n(M);
\end{equation}
\item
For any aggregation procedure $\bar f_n$, there exists a dictionary
$F$ of cardinality $M$ and a random variable $Z$ such that
$|Q(Z,f)| \leq b$ almost surely for all $f\in F$ and
\[
\E R(\bar f_n)\geq\min_{f\in F}R(f)+c_1\psi_n(M).
\]
\end{enumerate}
\end{Definition}
In our setup, one can show (cf.~\cite{TsyCOLT07}) that in
general, an optimal rate of aggregation (in the sense of
\cite{TsyCOLT07} [optimality in expectation] and of~\cite{LM1}
[optimality in probability]) is lower-bounded by $(\log M)/n$.
Thus, procedures satisfying an exact oracle inequality like
(\ref{eq:exact-oracle-inequality})---that is, an oracle inequality
with a factor of 1 in front of $\min_{f\in F}R(f)$---with a
residual term of $\psi_n(M)=(\log M)/n$ are said to be optimal.
Only a few aggregation procedures have been shown to achieve this
optimal rate, including the exponential aggregating schemes of \cite
{MR2163920,PhD04a,MR1762904,Audibert1,JRT08}, the the ``empirical
star algorithm'' in~\cite{Audibert1}, and the
``preselection/convexification algorithm'' in~\cite{LM1}. For a
survey on optimal aggregation procedures, see
the HDR dissertation of J.-Y. Audibert.

Our main focus here is on the problem of the optimality of the
aggregation procedure with exponential weights (AEW). This procedure
originate from the thermodynamic standpoint of
learning theory (see~\cite{MR2483528} for the state of the art in
this direction). AEW can be viewed as a relaxed version of the
trivial aggregation scheme, which is to minimize the empirical
risk
%
%
\begin{equation} \label{eq:R-n}
R_n(f)=\frac{1}{n}\sum_{i=1}^n Q(Z_i,f)
\end{equation}
in the dictionary $F$.

A procedure that minimizes \eqref{eq:R-n} is called
\textit{empirical risk minimization} (ERM). It is well known
that ERM generally cannot achieve the optimal rate of $(\log
M)/n$, unless one assumes that the given class $F$ has certain
geometric properties, which we discuss below (see also
\cite{LM2,m:08,JRT08}). To have any chance of obtaining better
rates, one has to consider aggregation procedures that take
values in larger sets than $F$. The most natural set is the
convex hull of $F$. AEW is a very popular candidate for the optimal
procedure, and it was one of the first procedures to be
studied in the context of the aggregation framework
\cite{JRT08,AoS08,LecAoS,MR2242356,MR2163920,PhD04a,MR1762904,DT07}.
It is defined by the following convex sum:
%
%
\begin{equation}
\label{eq:AEW}
\tilde{f}^{\mathrm{AEW}}=\sum_{j=1}^M\hat\theta_j f_j,\qquad \mbox{where }
\hat\theta_j=\frac{\exp(-({n}/{T})R_n(f_j))}
{\sum_{k=1}^M\exp(-({n}/{T})R_n(f_k))}
\end{equation}
for the dictionary $F=\{f_1,\ldots,f_M\}$.
The parameter $T>0$ is called the
\textit{temperature}.\footnote{This terminology comes
from thermodynamics, since the weights
$(\hat\theta_1,\ldots,\hat\theta_M)$ can be seen as a Gibbs
measure with temperature $T$ on the dictionary $F$.}

Thus far, there have been three main results concerning the
optimality of the AEW. The first of these is that the progressive
mixture rule
is optimal in expectation for $T$ larger than some parameters of
the model (see~\cite{MR2163920,MR1790617,MR1946426,JRT08,AoS08} and
\cite{Audibert1}),
and under certain convexity assumption on the loss function $Q$.
This procedure is defined by
%
%
\begin{equation}\label{eq:progressive-mixture}
\bar f=\frac{1}{n}\sum_{k=1}^n\tilde f^{\mathrm{AEW}}_k,
\end{equation}
where $\tilde f^{\mathrm{AEW}}_k$ is the function generated by AEW
(with a
common temperature parameter~$T$) associated with the dictionary
$F$ and constructed using only the first $k$ observations
$Z_1,\ldots,Z_k$. (See~\cite{Audibert1} for more details and for
other procedures related to the progressive mixture rule.)\vadjust{\goodbreak}

Second, the optimality in expectation of AEW was obtained by \cite
{DT07} for the regression model $Y_i=f(x_i)+\varepsilon_i$ with
a deterministic design $x_1,\ldots,x_n\in\cX$ with respect to the
risk $\norm{g-f}_n^2=n^{-1}\sum_{i=1}^n(g(x_i)-f(x_i))^2$ (with
its empirical version being
$R_n(g)=n^{-1}\sum_{i=1}^n(Y_i-g(x_i))^2$). That is, it was shown
that for $T\geq c\max(b,\sigma^2)$, where $\sigma^2$ is the
variance of the noise $\varepsilon$,
%
%
\begin{equation}\label{eq:Sacha-arnak}
\E\norm{\tilde f^{\mathrm{AEW}}-f}_n^2\leq\min_{g\in F}\norm
{g-f}_n^2+\frac{T\log M}{n+1}.
\end{equation}
Finally,~\cite{PhDPierre,PhD04a}, and~\cite{MR2483528}
proved that in the high-temperature regime, AEW can
achieve the optimal rate $(\log M)/n$ under the Bernstein
assumption, recalled below in Definition~\ref{def:Bernstein} in
expectation and in high probability. This result is discussedin more
detail later.

Despite the long history of AEW, the literature contains no results on the
optimality (or suboptimality) of AEW in the regression model with
random design in the general case (when the dictionary does not
necessarily satisfy the Bernstein condition). In this article, we
address this issue and complement the results (assuming the
Bernstein condition) of~\cite{PhDPierre,PhD04a,MR2483528} for the
low-temperature regime by proving the following:
\begin{enumerate}
\item[-]AEW is suboptimal for low temperatures $T\leq c_1$ (where $c_1$
is an absolute positive constant), both in expectation and in
probability, for the quadratic loss function and a dictionary of
cardinality $2$ (Theorem~\ref{thA}).
\item[-] AEW is suboptimal in probability for some large dictionaries
(of cardinality $M\sim\sqrt{n\log n}$) and small temperatures $T\leq
c_1$ (Theorem~\ref{thB}).
\item[-] AEW achieves the optimal rate $(\log M)/n$ for low
temperatures under the Bernstein condition on the dictionary (Theorem
\ref{thC}). Together with the high-temperature results of \cite
{PhDPierre,PhD04a} and~\cite{MR2483528}, this proves that the temperature
parameter has almost no impact (as long as $T=\cO(1)$) on the
performance of the AEW under this condition, with a residual term of
the order of $((T+1)\log M)/n$ for every $T>0$.
\end{enumerate}

\renewcommand{\thetheorema}{\Alph{theorema}}
\begin{theorema}\label{thA}There exist absolute constants
$c_0,\ldots,c_5$ for which the following holds. For any integer
$n\geq c_0$, there are random variables $(X,Y)$ and a dictionary
$F=\{f_1,f_2\}$ such that $(Y-f_i(X))^2\leq1$ almost surely for
$i=1,2$, for which the quadratic risk of the AEW satisfies the following:
\begin{enumerate}
\item if $T\leq c_1$ and $n$ is odd, then
\[
\E R(\tilde{f}^{\mathrm{AEW}})\geq\min_{f\in F}R(f)+\frac
{c_2}{\sqrt{n}};
\]
\item if $T \leq c_3\sqrt{n}/\log n$, then, with probability greater
than $c_4$,
\[
R(\tilde{f}^{\mathrm{AEW}})\geq\min_{f\in F}R(f)+\frac{c_5}{\sqrt{n}}.
\]
\end{enumerate}
\end{theorema}

Theorem~\ref{thA} proves that AEW is suboptimal in expectation in the
low-temperature regime and suboptimal in probability in both the\vadjust{\goodbreak} low- and
high-temperature regimes, since it is possible to construct
procedures that achieve the rate $C/n$ with high probability \cite
{Audibert1,LM1} and in expectation \cite
{MR2163920,MR1790617,MR1946426,JRT08,AoS08,Audibert1} in
the same setup as for Theorem~\ref{thA}. It should be noted that the problem
of the optimality in probability of the progressive mixture rule
(and other related procedures) was studied by~\cite{Audibert1},
who proved that, for a loss function $Q$ satisfying some
convexity and regularity assumption (e.g., the quadratic
loss used in Theorem~\ref{thA}), the progressive mixture rule $\bar f$
defined in (\ref{eq:progressive-mixture}) satisfies that for any
temperature parameter, with probability greater than an absolute
constant $c_0>0$, $R(\bar{f})\geq\min_{f\in F}R(f)+c_1 n^{-1/2}.$

In addition, it is important to observe that the suboptimality in
probability does not imply suboptimality in expectation for the
aggregation problem, or vice versa. This property of the
aggregation problem was first noted by~\cite{Audibert1}, who
found the progressive mixture rule (and other related aggregation
procedures) to be suboptimal in probability for
dictionaries of cardinality two but, on the other hand, to be
optimal in expectation (\cite{MR2163920,MR1790617,MR1946426} and
\cite
{JRT08}). This peculiar property of the problem
of aggregation comes from the fact that an aggregate $\hat f$ is
not restricted\vspace*{2pt} to the set $F$, which allows $R(\hat f
)-\min_{f\in F}R(f)$ to take negative values.~\cite{Audibert1} showed
that for the progressive mixture
rule $\bar f$, these negative values do compensate on average for
larger values, but there is still an event of constant probability
on which $R(\bar f )-\min_{f\in F}R(f)$ takes values greater than
$C/\sqrt{n}$.

The proof of Theorem~\ref{thA} shows that a dictionary consisting of two
functions is sufficient to yield a lower bound in expectation in the
low-temperature regime and in probability in both the
small temperature regime, $0\leq T\leq c_1$, and the large temperature
regime, $c_1\leq T\leq c_3\sqrt{n}/\log n$. In the following theorem,
we study the behavior of AEW for larger dictionaries. To the best of our
knowledge, negative results on the behavior of exponential weights
based aggregation procedures are not known for dictionaries with
more than two functions, and we show that the behavior of
the AEW deteriorates in some sense as the cardinality of the
dictionary increases.

\begin{theorema}\label{thB} There exist an integer $n_0$ and
absolute constants $c_1$ and $c_2$ for which the following holds.
For every $n \geq n_0$, there are random variables $(X,Y)$ and a
dictionary $F=\{f_1,\ldots,f_M\}$ of cardinality, $ M=\lceil
c_1\sqrt{n\log n}\rceil$, for which the quadratic loss function of
any element in $F$ is bounded by $2$ almost surely, and for every
$0<\alpha\leq1/2$, if $T \leq c_2\alpha$, then with probability
at least $1-c_3(\alpha)n^{\alpha-1/2}$,
\[
R(\tilde{f}^{\mathrm{AEW}})\geq\min_{f \in F} R(f)+c_4(\alpha
)\sqrt
{\frac{\log M}{n}}.
\]
Moreover, if $f_{F}^*\in F$ denotes the optimal function in $F$ with
respect to the quadratic loss (the oracle), then there exists $f_j \not
= f_F^*$ with an excess risk greater than $c_5(\alpha)n^{-1/2}$ and
for which
the weight of $f_j$ in the AEW procedure satisfies $\hat\theta_j \geq
1-n^{-c_6(\alpha)/T}.$
\end{theorema}

Theorem~\ref{thB} implies that the AEW procedure might cause the weights to
concentrate around a ``bad'' element in the dictionary (i.e., an
element whose risk is larger than the best in the class by at least
$\sim\!n^{-1/2}$) with high probability. In particular, Theorem~\ref{thB}
provides additional evidence that the AEW procedure is suboptimal for
low temperatures.

The analysis of the behavior of AEW for a dictionary of
cardinality larger than two is considerably harder than in the
two-function case and requires some results on rearrangement
of independent random variables that are almost Gaussian (see
Proposition~\ref{prop:Shahar-Proposition} below). Fortunately, not all
is lost as far as optimality results for AEW go. Indeed, we show that
under some geometric condition, AEW can be optimal and in fact can even
adapt to the ``real complexity'' of the dictionary.

Intuitively, a good aggregation scheme should be able to ignore
the elements in the dictionary whose risk is far from the optimal
risk in $F$, or at least the impact of such elements on the
function produced by the aggregation procedure should be small.
Thus, a good procedure is one with a residual term of the order
of $\psi/n$, where $\psi$ is a complexity measure that is
determined only by the richness of the set of ``almost minimizers''
in the dictionary. This leads to the following question:

\begin{Question}
\label{qu:adaptivityToComplexityOfDictionary} Is it possible to
construct an aggregation procedure that adapts to the real
complexity of the dictionary?
\end{Question}

This question was first addressed by the PAC-Bayesian approach. \cite
{PhDPierre,PhD04a} and~\cite{MR2483528}
showed that in the high-temperature regime, AEW satisfies the
requirements of
Question~\ref{qu:adaptivityToComplexityOfDictionary}, assuming
that the class has a geometric property, called the Bernstein
condition.

\begin{Definition}[(\cite{MR2240689})]\label{def:Bernstein}
We say that a function class $F$ is a $(\beta,B)$-Bernstein class
($0<\beta\leq1$ and $B\geq1$) with respect to $Z$ if every $f \in
F$ satisfies $\E f\geq0$ and
%
%
\begin{equation}
\label{eq:MarginAssumption}
\E(f^2(Z))\leq B(\E f(Z))^{\beta}.
\end{equation}
\end{Definition}

There are many natural situations in which the Bernstein condition
is satisfied. For instance, when $Q$ is the quadratic loss
function and the regression function is assumed to belong to $F$,
the excess loss function class
$\cL_F=\{Q(\cdot,f)-Q(\cdot,f_F^*)\dvt f\in F\}$ satisfies the
Bernstein condition with $\beta=1$, where $f_F^*\in F$ is the
minimizer of the risk in the class~$F$. Another generic example
is when the target function $Y$ is far from the set of targets
with ``multiple minimizers'' in $F$ and
$\cL_F$ satisfies the Bernstein condition with $\beta=1$. (See
\cite{m:08,MR2426759} for an exact formulation of this statement
and related results.)

The Bernstein condition is very natural in the context of ERM
because it has two consequences: that the empirical
excess risk has better concentration properties around the excess
risk, and that the complexity of the subset of $F$ consisting
of almost minimizers is smaller under this assumption. Consequently, if
the class $\cL_F$ is a $(\beta,B)$-Bernstein
class for $0 < \beta\leq1$, then the ERM algorithm can achieve
fast rates (see, e.g.,~\cite{MR2240689} and references
therein). As the results below show, the same is true for AEW.
Indeed, under a Bernstein assumption,
\cite{PhDPierre,PhD04a} and~\cite{MR2483528} proved that if
$R(\cdot)$ is a convex risk function and if $F$ is such that
$|Q(Z,f)|\leq b$ almost surely for any $f\in F$, then for every $T
\geq c_1\max\{b,B\}$ and $x>0$, with probability greater than
$1-2\exp(-x)$,
%
%
\begin{equation}\label{eq:PAC-Bound}
R(\tilde{f}^{\mathrm{AEW}})\leq\min_{f\in F}R(f)+\frac
{Tc_2}{n}\biggl(x+\log
\biggl(\sum_{f\in F}\exp\bigl(-(n/2T)\bigl(R(f)-R(f^*_F)\bigr)\bigr)\biggr)\biggr).
\end{equation}

Although the PAC-Bayesian approach cannot be used to obtain
\eqref{eq:PAC-Bound} in the low-temperature regime ($T \leq
c_1\max\{b,B\}$), such a result is not surprising. Indeed, because
fast error rates for the ERM are expected when the
underlying excess loss functions class satisfies the Bernstein
condition, and because AEW converges to the ERM when the temperature
$T$ tends to 0, it is likely that for ``small values'' of $T$,
AEW inherits some of the properties of ERM, such as fast
rates under a Bernstein condition. We show this in Theorem~\ref{thC}, proving that AEW answers
Question~\ref{qu:adaptivityToComplexityOfDictionary} for low temperatures
under the Bernstein condition.

Before formulating Theorem~\ref{thC}, we introduce the following
measure of complexity. For every $r>0$, let
\begin{eqnarray*}
\psi(r) &= & \log\bigl(|\{ f \in F \dvt R(f) - R(f_F^*) \leq r \}|+1\bigr)
\\
&&{}+  \sum_{j=1}^\infty2^{-j} \log\bigl(|\{ f \in F \dvt 2^{j-1} r < R(f) -
R(f_F^*) \leq2^j r \}|+1\bigr),
\end{eqnarray*}
where $|A|$ denotes the cardinality of the set $A$.

Observe that $\psi(r)$ is a weighted sum of the number of elements in $F$
that assigns smaller and smaller weights to functions with a relatively
large excess risk.

\begin{theorema}\label{thC} There exist absolute constants
$c_0$, $c_1,c_2$, and $c_3$ for which the following holds. Let~$F$ be
a class of functions bounded by $b$ such that the excess loss class
$\cL_F$ is a $(1,B)$-Bernstein class with respect to $Z$. If the risk
function $R(\cdot)$ is convex and if $T \leq c_0\max\{b,B\}$, then for
every $x> 0$, with probability at least $1-2\exp(-x)$, the function
$\tilde{f}^{\mathrm{AEW}}$ produced by the AEW algorithm satisfies
\[
R(\tilde{f}^{\mathrm{AEW}})\leq R(f^*_F) + c_1(b+B)\frac{x+ \psi
(\theta)}{n},
\]
where $\theta=c_2(b+B)(\log|F|)/n$.

In particular,
\[
\E R(\tilde{f}^{\mathrm{AEW}})\leq R(f^*_F) + c_3 (b+B)\frac{\psi
(\theta
)}{n} .
\]
\end{theorema}

In other words, the scaling factor $\theta$ that we use is proportional
to $(b+B)(\log|F|)/n$, and if the class is regular (in the sense
that the complexity of $F$ is well spread and not concentrated
just around one point), then $\psi(\theta)$ is roughly the cardinality
of the elements in $F$ with risk at most $\sim\!(b+B)(\log
|F|)/n$.

Observe that for every $r>0$, $\psi(r) \leq c \log|F|$ for a
suitable absolute constant $c$. Thus, if $T$ is reasonably
small (below a level proportional to $\max\{B,b\}$), then the
resulting aggregation rate is the optimal one, proportional to $
(b+B)(x+ \log M)/n$ with probability $1-2\exp(-x)$,
and proportional to $(b+B)(\log M)/n$ in expectation.
Thus, Theorem~\ref{thC} indeed gives a positive answer to Question
\ref{qu:adaptivityToComplexityOfDictionary} in the presence of a
Bernstein condition and for low temperatures.

Although the residual terms in Theorem~\ref{thC} and in
\eqref{eq:PAC-Bound} are not the same, they are comparable.
Indeed, the contribution of each element in $F$ in the residual
term depends exponentially on its excess risk.

Theorem~\ref{thC} together with the results for high temperatures from
\cite{PhDPierre,PhD04a} and~\cite{MR2483528} show that
the AEW is an optimal aggregation procedure under the Bernstein
condition as long as $T=\cO(1)$ when $M$ and $n$ tend to infinity.
In general, the residual term obtained is on the order of $
((T+1)\log M)/n$, and it can be proven that the optimal rate of
aggregation under the Bernstein condition is proportional to
$(\log M)/n$ using the classical tools in~\cite{MR2724359}.

Finally, a word about the organization of the article. In the next
section we present some comments about our results. The proofs of
the three theorems follow in the subsequent sections. Throughout, we
denote absolute constants or constants that depend on other
parameters by $c_1$, $c_2$, etc. (Of course, we specify
when a constant is absolute and when it depends on other
parameters.) The values of constants may change from line to line.
We write $a \sim b$ if there are absolute constants $c$ and $C$
such that $bc \leq a \leq Cb$, and write $a \lesssim b$ if $a \leq Cb$.

\section{Comments}
\label{sec:Comments}

Although from a theoretical standpoint, whether AEW is an optimal
procedure in expectation and for high
temperatures in the regression model with random design remains to be
seen, from a
practical standpoint, we believe that exponential aggregating
schemes simply should not be used in the setup of this article,
because of the following reasons (see also the comments in~\cite{Audibert1}):
\begin{enumerate}
\item For any temperature $T\leq c_0\sqrt{n}/\log n$, there is an event
of constant probability on which AEW performs poorly (this is the
second part of Theorem~\ref{thA}).
\item If the temperature parameter is chosen to be too small, then the
AEW can perform poorly even in expectation (the first part of Theorem~\ref{thA}).
\end{enumerate}

Another consequence of the lower bounds stated in Theorem~\ref{thA} is
that AEW cannot be an optimal aggregation procedure both in
expectation and in probability at low temperatures for two other
aggregation problems: the problem of \textit{convex aggregation},
in which one wants to mimic the best element in the convex hull of
$F$, and the problem of \textit{linear aggregation}, where one
wishes to mimic the best linear combination of elements in $F$.
Indeed, clearly
\[
\min_{f\in F}R(f)\geq\min_{f\in\operatorname{conv}(F)}R(f)\geq
\min
_{f\in\operatorname{span}(F)}R(f).
\]
Moreover, the optimal rates of aggregation for the convex and linear
aggregation problems for dictionaries of cardinality two are of
the order of $n^{-1}$ (see~\cite{TsyCOLT07,MR2329442,LM7}), whereas
the residual terms obtained in Theorem~\ref{thA} are on the order of
$n^{-1/2}$ for such a dictionary. Thus AEW is suboptimal for
these two other aggregation problems in the low-temperature
regime.

We end this section by comparing two seemingly related
assumptions, the margin assumption of~\cite{Tsy04} and the
Bernstein condition of~\cite{MR2240689}. Note that in
the proof of Theorem~\ref{thC}, we have restricted ourselves to the case
$\beta=1$ simply to make the presentation as simple as possible.
A~very similar result, with the residual term
$((x+\psi(\theta))/n)^{1/(2-\beta)}$ for the exact oracle
inequality in probability and
$(\psi(\theta)/n)^{1/(2-\beta)}$ for the exact oracle
inequality in expectation, holds if one assumes a Bernstein
condition for any $0<\beta<1$, and the proof is identical to that in
the case where $\beta=1$. This makes the discussion about
$\beta$-Bernstein classes relevant here.

Recall the definition of the margin assumption:
\begin{Definition}[(\cite{Tsy04})]
We say that $F$ has margin with parameters $(\beta,B)$ ($0<\beta\leq1$
and $B\geq1$) if for every $f \in F$,
\[
\E\bigl(\bigl(Q(Z,f)-Q(Z,f^*)\bigr)^2\bigr)\leq B\bigl(R(f)-R(f^*)\bigr)^\beta,
\]
where $f^*$ is defined such that $R(f^*)=\min_f R(f)$, and the
minimum is taken with respect to all measurable functions $f$ on the
given probability space.
\end{Definition}

Although the margin condition appears similar to the Bernstein
condition, they are in fact very different, and have been
introduced in the context of different types of problems.
In the first of these, the ``classical'' statistical setup, one is
given a
function class $F$ (the \textit{model}) with an upper bound on its
complexity and an unknown target function $f^*$, the
minimizer of the risk over \textit{all} measurable functions. One
usually assumes that $f^*$ belongs to $F$, and the aim is to
construct an estimator $\hat{f}=\hat f (\cdot,\cD)$ for which the
risk $R(\hat f)$ tends to 0 quickly as the sample size tends to
infinity. In this setup, the margin assumption can improve this
rate of convergence because of a better concentration of empirical
means of $Q(\cdot,f)-Q(\cdot, f^*)$ around its mean~\cite{Tsy04}.
The margin assumption (MA) for $\beta=1$ compares the
performance of each $f \in F$ with the \textit{best possible measurable
function}, but it has nothing to do with the geometric structure
of $F$. The margin is determined for every $f$ separately, because
$f^*$ does not depend on the choice of $F$.

In the second type of problem, the ``learning theory'' setup, one does
not assume
that the target function $f^*$ belongs to $F$. The aim is to
construct a function $\hat f$ with a risk as close as possible
to that of the best element $f_F^* \in F$. Assuming that the
excess loss class $\cL_F$ satisfies the Bernstein condition (BC), the
error rate can be improved (see, e.g.,~\cite{MR2426759,MR2240689}).

At a first glance, MA and BC (for $\beta=1$) share very strong
similarities. Indeed, saying that $\cL_F$ is a $(1,B)$-Bernstein
class means that for every $f\in F$,
\[
\E\bigl(\bigl(Q(Z,f)-Q(Z,f_F^*)\bigr)^2\bigr)\leq B\bigl(R(f)-R(f_F^*)\bigr),
\]
but nevertheless they are different. Indeed, as mentioned earlier, MA
is only a matter of concentration (and classical statistics
questions are mostly a question of the trade-off between
concentration and complexity). On the other hand, BC involves a
lot of geometry of the function class $F$, because $f^*_F$ might
change significantly by adding a single function to $F$ or by
removing a function. In fact, the difficulty of learning theory
problems is determined by the trade-off between concentration and
complexity, \textit{and} the geometry of the given class, since one
measures the performance of the learning algorithm relative to the
best \textit{in the class}. Assuming that $f^* \in F$, as is usually
done in classical statistics, exempts one from the need to
consider the geometry of~$F$, but one does not have that freedom in
the aggregation framework. Indeed, since in the AEW algorithm the
estimator is determined by the empirical means
$R_n(f)-R_n(f^*_F)$, this is a learning problem rather than a
problem in classical statistics, despite the fact that it has been
used in statistical frameworks to construct adaptive estimators (see,
e.g.,~\cite{AoS08,GL1,LecAoS,ST07,BTW07,MR2242356,Tsy04,PhD04a,MR1762904}).
Therefore, given their nature, aggregation procedures like the
AEW are more natural under a BC assumption than under the MA.
(A by-product of Theorem~\ref{thA} is that the MA cannot improve the
performance of AEW since in the setup of Theorem~\ref{thA}, it is easy to check
that MA is satisfied with the best possible margin parameter
$\beta=1$.)

\section{Preliminary results on Gaussian approximation}
\label{sec:preliminaries} Our starting point is the
Berry--Ess{\'e}en theorem on Gaussian approximation. Let
$(W_n)_{n\in\N}$ be a sequence of i.i.d., mean-0 random
variables with variance $1$, set $g$ to be a standard Gaussian
variable, and write
\[
\bar X_n=\frac{1}{\sqrt{n}}\sum_{i=1}^nW_i.
\]

\begin{Theorem}[(\cite{MR1353441})]\label{theo:Berry-Esseen}There exists
an absolute constant $A>0$ such that for every integer~$n$,
\[
\sup_{x\in\R}|\Pro[\bar X_n\leq x]-\Pro[g\leq x]|\leq\frac{A\E
|W_1|^3}{\sqrt{n}}.
\]
\end{Theorem}
From here on, we let $A$ denote the constant appearing in Theorem \ref
{theo:Berry-Esseen}.

When the tail behavior of the $W_i$ has a subexponential decay, the
Gaussian approximation can be improved. Indeed, recall that a
real-valued random variable $W$ belongs to $L_{\psi_\alpha}$ for some
$\alpha\geq1$ if there exists $0<c<\infty$ such that
%
%
\begin{equation}
\label{eq:psi-alpha}
\E\exp(|W|^\alpha/c^\alpha)\leq2.
\end{equation}
The infimum over all constants $c$ for which (\ref{eq:psi-alpha}) holds
defines an Orlicz norm, which is called the $\psi_\alpha$ norm and is
denoted by $\| \cdot\|_{\psi_\alpha}$. (For more information on Orlicz
norms, see, e.g.,~\cite{vanderVaartWellner} and~\cite{MR1113700}.)

\begin{Proposition}[(Chapter 5 in~\cite{MR1353441})]\label{prop:Petrov}
For every $L>0$, there exist constants $B_0, c_1$, and $c_2$ that
depend only on $L$ for which the following holds. If $\|W\|_{\psi_1}
\leq L$, then for any $x\geq0$, such that $x \leq B_0n^{1/6}$,
\[
\Pro[\bar X_n\geq x]=\Pro[g\geq x]\exp\biggl(\frac{x^3 \E W^3 }{6\sqrt
{n}}\biggr)\biggl[1+\mathrm{O}\biggl(\frac{x+1}{\sqrt{n}}\biggr)\biggr]
\]
and
\[
\Pro[\bar X_n\leq-x]=\Pro[g\leq-x]\exp\biggl(-\frac{x^3\E W^3 }{6\sqrt
{n}}\biggr)\biggl[1+\mathrm{O}\biggl(\frac{x+1}{\sqrt{n}}\biggr)\biggr],
\]
where by $v=\mathrm{O}(u)$ we mean that $-c_1u \leq v \leq c_1u$.

In particular, if $|x| \leq B_0n^{1/6}$ and $\E W^3=0$, then
\[
|\Pro[\bar X_n\leq x]-\Pro[g\leq x]| \leq c_2\bigl(n^{-1/2}\exp(-x^2/2)\bigr).
\]
\end{Proposition}
From here on, we let $B_0$ denote the constant appearing in
Proposition~\ref{prop:Petrov}.

\section{\texorpdfstring{Proof of Theorem \protect\ref{thA}}{Proof of Theorem A}}\label{sec:proof-of-theo-A}
Before presenting the proof of Theorem~\ref{thA}, we introduce the following notation. Given a probability
measure $\nu$ and $(Z_i)_{i=1}^n$ selected independently according
to $\nu$, we set $P_n = n^{-1}\sum_{i=1}^n \delta_{Z_i}$ the
empirical measure supported on $(Z_i)_{i=1}^n$. We let $P$
denote the expectation $\E_\nu$. We assume that $T \leq1$ and
recall that $n$ is an odd integer.

Let $Y=0$ and define $X$ by $\Pro[X=1]=1/2-n^{-1/2}$ and $\Pro
[X=-1]=1/2+n^{-1/2}$. Let $f_1=\1_{[0,1]}$ and $f_2=\1_{[-1,0]}$, and
consider the dictionary $F=\{f_1,f_2\}$. It is easy to verify that the
best function in $F$ (the oracle) with respect to the quadratic risk is
$f_1$, and that the excess loss function of $f_2$, $\cL
_2=f_2^2-f_1^2=f_2-f_1$, satisfies that
\[
\cL_2(X)=-X,\qquad \E\cL_2(X)=2n^{-1/2}\quad \mbox{and}\quad \sigma^2=\E(\cL
_2(X)-\E\cL
_2(X))^2=1-4/n.
\]
To simplify notation, set $P\cL_2=\E\cL_2(X)$ and
$P_n\cL_2=n^{-1}\sum_{i=1}^n\cL_2(X_i)$.

An important parameter that lies at the heart of this
counterexample is the Bernstein constant (which is very bad in
this case),
%
%
\begin{equation}\label{eq:def-alpha}
\alpha=\frac{\E(f_1-f_2)^2}{P \cL_2}=\frac{\sqrt{n}}{2}.
\end{equation}
Straightforward computation shows that
AEW on $F$ with temperature $T$ is given by
\[
\tilde{f}^{\mathrm{AEW}}=\hat\theta_1f_1+(1-\hat\theta_1)f_2,\qquad
\hat\theta
_1=\frac{1}{1+\exp(-({n}/{T})P_n\cL_{2})},
\]
and that for $h(\theta)=\theta+\alpha\theta(1-\theta)$ defined
for all $\theta\in[0,1]$,
%
\begin{eqnarray}
\label{eq:excess-risk-AEW}
\E[R(\tilde{f}^{\mathrm{AEW}})-R(f_1)]&=&\E[1-\hat\theta
_1-\alpha\hat\theta_1(1-\hat\theta_1)]P\cL_2=\E[1-h(\hat\theta
_1)]P\cL_2
\nonumber\\
&=&\biggl[1-\int_0^\infty h^\prime(t)\Pro[\hat\theta_1\geq t]\,
\mathrm{d}t\biggr]P\cL_2
\nonumber
\\[-8pt]
\\[-8pt]
\nonumber
&=&\biggl[1+\int_0^1\bigl(2\alpha t-(1+\alpha)\bigr)\Pro[\hat
\theta
_1\geq t]\,\mathrm{d}t\biggr]P\cL_2
\\
&=&\biggl[1+\int_0^1\bigl(2\alpha t-(1+\alpha)\bigr)\Pro[ P_n\cL_2\geq\gamma(t) ]\,
\mathrm{d}t\biggr]P\cL_2,\nonumber
\end{eqnarray}
where $\gamma(t)$ is an increasing function defined for any $t\in
(0,1)$ by
\[
\gamma(t)=\frac{T}{n}\log\biggl(\frac{t}{1-t}\biggr).
\]
In particular,
\[
\E[R(\tilde{f}^{\mathrm{AEW}})-R(f_1)]=[I_1+I_2]P\cL_2
\]
for
\[
\label{eq:I1}
I_1=\int_0^{\alpha^{-1}}\bigl(2\alpha t-(1+\alpha)\bigr)\Pro[P_n\cL_2\geq
\gamma
(t)]\,\mathrm{d}t+1
\]
and
\[
\label{eq:I2}
I_2=\int_{\alpha^{-1}}^{1}\bigl(2\alpha t-(1+\alpha)\bigr)\Pro[P_n\cL_2\geq
\gamma
(t)]\,\mathrm{d}t.
\]

First, we bound $I_1$ from below. To that end, we note the following
facts. First, for every $0\leq t\leq\alpha^{-1}$, $1+\alpha-2\alpha
t\geq0$ and
\[
\int_0^{\alpha^{-1}}\bigl(2\alpha t -(1+\alpha)\bigr)\,\mathrm{d}t=-1.
\]
Second, if we set $E=\exp(n P\cL_2/T)$, then for $T \lesssim\sqrt
{n}/\log n$, $0<(1+E)^{-1}\leq\alpha^{-1}$. In particular, this holds
under our assumption that $T\leq1$. Moreover, because $\gamma$ is
increasing, for $(1+E)^{-1}\leq t\leq\alpha^{-1}$, $\gamma(t) \geq
\gamma((1+E)^{-1})= -P\cL_2$. Therefore,
\begin{eqnarray*}
I_1&=&\int_0^{\alpha^{-1}}\bigl(2\alpha t-(1+\alpha)\bigr)\Pro[P_n\cL_2\geq
\gamma
(t)]\,\mathrm{d}t+1
\\
&=&\int_0^{\alpha^{-1}}\bigl(2\alpha t-(1+\alpha)\bigr)\bigl(\Pro[P_n\cL_2\geq
\gamma
(t)]-1\bigr)\,\mathrm{d}t
\\
&\geq&\int_{(1+E)^{-1}}^{\alpha^{-1}}(1+\alpha-2\alpha t)\Pro
[P_n\cL
_2<\gamma(t)]\,\mathrm{d}t
\\
&\geq&\int_{(1+E)^{-1}}^{\alpha^{-1}}(1+\alpha-2\alpha t)\,\mathrm{d}t
\cdot\Pro\bigl[\bigl(\sqrt{n}/\sigma\bigr)(P_n\cL_2-P\cL_2)<\bigl(\sqrt{n}/\sigma
\bigr)(-2P\cL_2)\bigr]
\\
&\geq&\int_{(1+E)^{-1}}^{\alpha^{-1}}(1+\alpha-2\alpha t)\,\mathrm
{d}t\bigl(\Pro[g\leq-8]-A/\sqrt{n}\bigr)\geq c_0>0,
\end{eqnarray*}
where in the last step we used the Berry--Ess{\'e}en theorem, with
$|\cL
_2| \leq1$ and $n\geq8\vee(2A/\Pro[g\leq-8])^2$, implying that $0<c_0<1/2$.

We turn to a lower bound for $I_2$. Applying a change of variables
$t\mapsto1+\alpha^{-1}-u$ in the second term of $I_2$, it is evident that
\begin{eqnarray*}
I_2&=&\int_{\alpha^{-1}}^{({\alpha+1})/({2\alpha})}\bigl(2\alpha
t-(1+\alpha
)\bigr)\Pro[P_n\cL_2\geq\gamma(t)]\,\mathrm{d}t
\\
&&{}+\int_{({\alpha+1})/({2\alpha})}^1\bigl(2\alpha t-(1+\alpha)\bigr)\Pro[P_n\cL
_2\geq
\gamma(t)]\,\mathrm{d}t
\\
&=&\int_{\alpha^{-1}}^{({\alpha+1})/({2\alpha})}\bigl(2\alpha t-(1+\alpha
)\bigr)\Pro
[\gamma(t)\leq P_n\cL_2<\gamma(1+\alpha^{-1}-t)]\,\mathrm{d}t=I_3+I_4
\end{eqnarray*}
for
\[
I_3= \int_{\alpha^{-1}}^{(1+c_0/4)\alpha^{-1}}\bigl(2\alpha t-(1+\alpha
)\bigr)\Pro
[\gamma(t)\leq P_n\cL_2<\gamma(1+\alpha^{-1}-t)]\,\mathrm{d}t
\]
and
\[
I_4=\int_{(1+c_0/4)\alpha^{-1}}^{({\alpha+1})/({2\alpha})}\bigl(2\alpha
t-(1+\alpha)\bigr)\Pro[\gamma(t)\leq P_n\cL_2<\gamma(1+\alpha
^{-1}-t)]\,
\mathrm{d}t.
\]
To estimate $I_3$, note that $2\alpha t -(1+\alpha) \leq0$ for $t \in
[\alpha^{-1},(\alpha+1)/(2\alpha)]$, and thus
\[
I_3 \geq\int_{\alpha^{-1}}^{(1+c_0/4)\alpha^{-1}} \bigl(2\alpha t
-(1+\alpha)\bigr)\,\mathrm{d}t \geq
\frac{-c_0}{4}\biggl(1+\frac{1}{\alpha}\biggr) \geq-\frac{c_0}{3}
\]
for our choice of $\alpha$.

The final step of the proof is to bound $I_4$ and in particular to show
that for small values of $T$, $I_4 \geq-c_0/3$.

For any $0<t\leq(\alpha+1)/(2\alpha)$, consider the intervals
$I_T(t)=[n\gamma(t),n\gamma(1+\alpha^{-1}-t))$,
and set $N_T(t)=|\{I_T(t)\cap\Z\}|$, which is the number of integers
in $I_T(t)$.
Because $\cL_2(X)=-X$,
\[
\Pro[\gamma(t)\leq P_n\cL_2<\gamma(1+\alpha^{-1}-t)]
=\Pro\Biggl[\sum_{i=1}^n-X_i\in I_T(t)\Biggr]=\Pro_T(t).
\]
Recall that $X \in\{-1,1\}$, and thus $\Pro[\sum_{i}-X_i\in
I_T(t)]=\Pro[\sum_{i}-X_i\in I_T(t) \cap\Z]$. Because $n\gamma(t)$ is
increasing and non-negative for $t>1/2$, then if $1/2<t \leq(\alpha
+1)/(2\alpha)$, it follows that
$0<n\gamma(t)<n\gamma(1+1/\alpha-t)<1$, provided that $T \leq1$. Thus,
for such values of $t$, $N_T(t)=0$, implying that $\Pro_T(t)=0$. On the
other hand, if $t\leq1/2$, then $\{0\}\subset I_T(t)\cap\Z$. In
particular, if $N_T(t)=1$, then $I_T(t)\cap\Z=\{0\}$, and since $n$ is
odd, then $\Pro_T(t)=\Pro[\sum_{i=1}^n-X_i=0]=0$. Otherwise,
$N_T(t)\geq
2$, which implies that $N_T(t)\leq2\Delta_T(t)$, where $\Delta_T(t)$
is the length of $I_T(t)$, given by
\[
\Delta_T(t)=n\bigl(\gamma(1+\alpha^{-1}-t)-\gamma(t)\bigr)=T\log\biggl(\frac
{(1-t)(\alpha+1-\alpha t)}{t(\alpha t-1)}\biggr).
\]
Therefore, for every $t$ in our range,
\[
\Pro_T(t)\leq N_T(t)\max_{k\in I_T(t)}\Pro\Biggl[\sum_{i=1}^n-X_i=k\Biggr]\leq
2\Delta_T(t)\max_{k\in\Z}\Pro\Biggl[\sum_{i=1}^nX_i=k\Biggr].
\]
Since $2\alpha t-(1+\alpha)\leq0$ for every $0<t\leq(\alpha
+1)/(2\alpha
)$, it is evident that
\[
I_4 \geq2 T \max_{k\in\Z}\Pro\Biggl[\sum_{i=1}^nX_i=k\Biggr] \cdot\int
_{(1+c_0/4)\alpha^{-1}}^{({\alpha+1})/({2\alpha})}
\bigl(2\alpha t-(1+\alpha)\bigr)\log\biggl(\frac{(1-t)(\alpha+1-\alpha t)}{t(\alpha
t-1)}\biggr) \,\mathrm{d}t.
\]
It can be shown that $\max_{k\in\Z}\Pro[\sum_{i=1}^n X_i=k]$ is on the
order of $n^{-1/2}$ either by a direct computation or by the
Berry--Ess{\'e}en theorem. Moreover, for any $(1+c_0/4)\alpha^{-1}\leq
t\leq(\alpha+1)/(2\alpha)$, one has $\alpha t-1\geq
c_0(4+c_0)^{-1}\alpha t$, and thus,
\[
\log\biggl(\frac{(1-t)(\alpha+1-\alpha t)}{t(\alpha t-1)}\biggr)\leq\log\biggl(\frac
{2(4+c_0)}{c_0 t^2}\biggr).
\]
Therefore, combining the two observations with a change of variables
$u=Ct$ for $C=(c_0/(2(4+c_0)))^{1/2}$, it is evident that there are
absolute constants $c_1,c_2$ for which
\[
I_4\geq\frac{c_1 T}{\sqrt{n}}\int_{C(1+c_0/4)\alpha
^{-1}}^{({C(\alpha
+1)}/(2\alpha))}(1+\alpha-2\alpha u/C)(\log u)\,\mathrm{d}u
\geq-c_2\frac{T\alpha}{\sqrt{n}}.
\]
Thus, there is an absolute constant $c_3$ such that if $T\leq c_3$,
then $I_4\geq-c_0/3$, implying that
\[
\E[R(\tilde{f}^{\mathrm{AEW}})-R(f_1)]\geq\frac{c_0}{3\sqrt{n}},
\]
and proving the first part of Theorem~\ref{thA}.

To prove the second part of the theorem, note that by the Berry--Ess{\'
e}en theorem, for every $x\in\R$, with probability greater than $\Pro
[g\leq x]-2A/\sqrt{n}$,
\[
\frac{\sqrt{n}}{\sigma(\cL_2)}(P_n\cL_2-P\cL_2)\leq x.
\]
Thus, if $n$ is large enough to ensure that
$\Pro[g\leq-4]-2A/\sqrt{n}\geq\Pro[g\leq-4]/2=c_4$,
and taking $x=-4$, then with probability at least $c_4$, $P_n\cL_2\leq
-n^{-1/2}$. In that case, $\hat\theta_1\leq\exp(-\sqrt{n}/T)$, which
yields that
\[
R(\tilde{f}^{\mathrm{AEW}})-R(f_1)=
\bigl(1-\hat\theta_1-\alpha\hat\theta_1(1-\hat\theta_1)\bigr)\cdot P\cL_2
\geq P\cL_2/4=n^{-1/2}/2,
\]
provided that $T \lesssim\sqrt{n}/\log n$.

\section{\texorpdfstring{Proof of Theorem \protect\ref{thB}}{Proof of Theorem B}}
\label{sec:Proof-of-theorem-B}

The first step in the proof of Theorem~\ref{thB} involves a general statement
regarding a monotone rearrangement of independent random variables that
are close to being Gaussian. Let $W$ be a mean~0, variance 1 random
variable that is absolutely continuous with respect to the Lebesgue
measure. Further assume that $|W|$ has a finite third moment (in fact,
the random variables in which we are interested are bounded) and set
$\beta(W)=A \E|W|^3$, where $A$ is the constant appearing in the
Berry--Ess{\'e}en theorem (Theorem~\ref{theo:Berry-Esseen}). Let
$W_1,\ldots,W_n$ be independent random variables distributed as $W$ and
set $\bar X= n^{-1/2} \sum_{i=1}^n W_i$. Let $(\bar X_j)_{j=1}^\ell$ be
$\ell$ independent copies of $\bar X$, and put $\gamma_1=\gamma
_1(\ell)
\in\R$ to satisfy that
\[
\Pro\Bigl[\min_{1\leq j\leq\ell}\bar X_j\leq\gamma_1(\ell)\Bigr]=1-\frac{1}{n}.
\]
Note that such a $\gamma_1$ exists because $W$ has a density with
respect to the Lebesgue measure.

Throughout the proof of Theorem~\ref{thB}, we require the following
simple estimates on $\gamma_1$.

\begin{Lemma} \label{lemma:gamma-1}
There exist absolute
constants $c_0,\ldots,c_3$ for which the following hold:
\begin{enumerate}
\item If $\ell\geq c_0\log n$, then
\[
1-\frac{\log n}{\ell} \leq\Pro[\bar X > \gamma_1] \leq1-c_1\frac
{\log
n}{\ell}.
\]
\item If $\ell$ and $n$ are such that $({\beta(W)}/{\sqrt{n}} +
({\log n})/{\ell} ) < \Pro[g < -2]$, then $\gamma_1 \leq-2$.

\item If $\gamma_1 \leq-2$ and $c_0\log n \leq\ell\leq
c_2\beta
^{-1}(W)\sqrt{n}\log n$, then
\[
|\gamma_1| \sim\log^{1/2}\biggl(\frac{c_3\ell}{\log n}\biggr) \quad\mbox{and}\quad
\exp(-\gamma_1^2/2) \sim\frac{\log n}{\ell}\log^{1/2}\biggl(\frac{c_3
\ell
}{\log n}\biggr).
\]
\end{enumerate}
\end{Lemma}

Before we present the proof of Lemma~\ref{lemma:gamma-1}, recall
that for every $x\geq2$,
%
%
\begin{equation}\label{eq:Gaussian-tail-estimate}
\frac{3}{4\sqrt{2\uppi}}\frac{\exp(-x^2/2)}{x}\leq\Pro[g\geq
x]\leq\frac{1}{\sqrt{2\uppi}}\frac{\exp(-x^2/2)}{x}.
\end{equation}

\begin{pf*}{Proof of Lemma \protect\ref{lemma:gamma-1}}
To prove the first part, note that by independence and because $\exp
(-x)\geq1-x$,
%
%
\begin{equation}\label{eq:Def-gamma1-2}
\Pro[\bar X>\gamma_1]=\Pro\Bigl[\min_{1 \leq j \leq\ell}\bar
X_j>\gamma_1\Bigr]^{{1}/{\ell}}=\biggl(\frac{1}{n}\biggr)^{1/\ell}
\geq1-\frac{\log n}{\ell}.
\end{equation}
The reverse inequality follows in an identical fashion, because $\exp
(-x)\leq1-x/3$ if $0 \leq x \leq1$.

Turning to the second part, if $\gamma_1 > -2$, then
\[
1-\frac{1}{n} =\Pro\Bigl[\min_{1 \leq j \leq\ell} \bar X_j \leq
-\gamma_1\Bigr] \geq\Pro\Bigl[\min_{1 \leq j \leq\ell} \bar X_j \leq-2\Bigr]
= 1-(\Pro[\bar X >-2])^\ell,
\]
implying that $\Pro[\bar X \leq-2] \leq(\log n)/\ell$. On the
other hand, by the Berry--Ess{\'e}en theorem, $\Pro[\bar X \leq-2]
\geq
\Pro[g \leq-2] - \beta(W)/\sqrt{n}$, which is impossible under
the assumptions of (2).

Finally, to prove (3), we use the Berry--Ess{\'e}en theorem
combined with the lower and upper estimates on the Gaussian tail
(\ref{eq:Gaussian-tail-estimate}) and (\ref{eq:Def-gamma1-2}). Thus,
\[
\frac{3}{4\sqrt{2\uppi}}\frac{1}{|\gamma_1|}\exp\biggl(-\frac{|\gamma
_1|^2}{2}\biggr)\leq
\Pro[g< \gamma_1]\leq\Pro[\bar
X<\gamma_1]+\frac{\beta(W)}{\sqrt{n}}\leq
\frac{\beta(W)}{\sqrt{n}}+c_1\frac{\log n}{\ell},
\]
and
\[
\frac{1}{\sqrt{2\uppi}}\frac{1}{|\gamma_1|}\exp\biggl(-\frac{|\gamma
_1|^2}{2}\biggr)\geq
\frac{\log n}{\ell}-\frac{\beta(W)}{\sqrt{n}},
\]
from which both parts of the third claim follow.
\end{pf*}

\begin{Proposition}\label{prop:Shahar-Proposition}
There exist constants $c_1,c_2,c_3$, and $c_4$ that depend only on
$\|W\|_{\psi_2}$ for which the following holds. Let $2M^2\exp(-c_1
n^{1/3})<\delta\leq1$, and assume that $\E W^3=0$ and that
$\gamma_1=\gamma_1(M-1)\leq-2$. Then
\begin{eqnarray*}
&&\Pro [ \exists j\in\{2,\ldots,M\}\dvt \bar X_j\leq\gamma_1
{\mbox{ and for every }} k\in\{2,\ldots,M\} \setminus\{j\}, \bar
X_k-\bar X_j\geq\delta]
\\
&&\quad \geq
1-\frac{1}{n}-c_2\biggl(\frac{1}{\sqrt{n}}+\delta\biggr)(\log
n)^2\sqrt{\log M},
\end{eqnarray*}
provided that $c_3\log n\leq M\leq c_4\sqrt{n}(\log n)$.
\end{Proposition}

\begin{pf}For every $2 \leq j \leq M$, let
\[
\Omega_j=\bigl\{\bar X_j\leq\gamma_1 {\mbox{ and }} \bar X_k- \bar X_j\geq
\delta\mbox{ for every }
k\in\{2,\ldots,M\} \setminus\{j\}
\bigr\}.
\]
The events $\Omega_j$ for $2 \leq j \leq M$ are disjoint, and thus
\begin{eqnarray*}
&& \Pro[\exists j\in\{2,\ldots,M\}\dvt\bar X_j\leq\gamma_1 {\mbox{ and
}} \bar X_k-\bar X_j\geq\delta{\mbox{ for every }} k\in\{2,\ldots,M\}
\setminus\{j\}]\\
&&\quad=\Pro\Biggl[\bigcup_{j=2}^M\Omega_j\Biggr]=(M-1)\Pro[\Omega_2].
\end{eqnarray*}
Since the variables $(\bar X_j)_{j=2}^M$ are independent, we have
\[
\Pro[\Omega_2]=\int_{-\infty}^{\gamma_1}f_{\bar
X}(z)\biggl(\int_{z+\delta}^\infty f_{\bar
X}(t)\,\mathrm{d}\mu(t)\biggr)^{M-2}\,\mathrm{d}\mu(z),
\]
where $f_{\bar X}$ is a density function of $\bar X$ with respect to
the Lebesgue measure $\mu$.

On the other hand, for any $z\leq\gamma_1$, $\Pro[\bar X\geq z]>0$
because of (\ref{eq:Def-gamma1-2}). Thus, for every
$z\leq\gamma_1$,
%
%
\begin{equation}
\label{eq:Proba-Omega2-2}
\int_{z+\delta}^\infty f_{\bar X}(t)\,\mathrm{d}\mu(t)=\biggl(1-\frac
{\int
_{z}^{z+\delta} f_{\bar
X}(t)\,\mathrm{d}\mu(t)}{\int_{z}^\infty f_{\bar X}(t)\,\mathrm
{d}\mu
(t)}\biggr) \cdot\int_{z}^\infty f_{\bar
X}(t)\,\mathrm{d}\mu(t).
\end{equation}
Note that for every $0 \leq x \leq1$, $(1-x)^{M-2}\geq1-(M-2)x$,
and applied to \eqref{eq:Proba-Omega2-2},
\begin{eqnarray*}
\Pro[\Omega_2]&\geq&\int_{-\infty}^{\gamma_1}f_{\bar X}(z)\biggl(\int
_z^\infty
f_{\bar X}(t)\,\mathrm{d}\mu(t)\biggr)^{M-2}\,\mathrm{d}\mu(z)
\\[-2pt]
&&{}-(M-2)\int_{-\infty}^{\gamma_1} f_{\bar X}(z)\biggl(\int_{z}^{\infty
}f_{\bar
X}(t)\,\mathrm{d}\mu(t)\biggr)^{M-3}
\biggl(\int_{z}^{z+\delta}f_{\bar X}(t)\,\mathrm{d}\mu(t)\biggr)\,\mathrm
{d}\mu(z)
\\[-2pt]
&\geq& \Pro[\bar X_2\leq\gamma_1 \mbox{ and } \bar X_k\geq\bar X_2,
\mbox{for every}
k \geq3]-T_2
\\[-2pt]
&=&\frac{1}{M-1}\Pro\Bigl[\min_{2\leq j\leq M}\bar
X_j\leq\gamma_1\Bigr]-T_2,
\end{eqnarray*}
where
\[
T_2=(M-2)\int_{-\infty}^{\gamma_1} f_{\bar
X}(z)\biggl(\int_{z}^{z+\delta}f_{\bar X}(t)\,\mathrm{d}\mu(t)\biggr)\,\mathrm
{d}\mu(z).
\]
Recall the if $(W_i)$ are independent mean-0 random variables and
$(a_i)$ are real numbers, then $\|\sum a_i W_i \|_{\psi_2} \leq c
(\sum
a_i^2 \|W_i\|_{\psi_2}^2 )^{1/2}$, where $c$ is an absolute constant
\cite{vanderVaartWellner}. Thus, $\|{\bar X}\|_{\psi_2}\leq c\norm
{W}_{\psi_2}$, and for any $t< 0$,
%
\[
\int_{-\infty}^{t} f_{\bar
X}(z)\biggl(\int_{z}^{z+\delta}f_{\bar X}(t)\,\mathrm{d}\mu(t)\biggr)\,\mathrm
{d}\mu
(z)\leq
\Pro[\bar X\leq t]\leq2\exp(-t^2/c^2\|{W}\|^2_{\psi_2}).
\]
Let $t_0<0$ be such that
\[
2\exp(-t_0^2/c^2\|{W}\|^2_{\psi_2})=\frac{\delta\sqrt{\log
(M-1)}}{(M-1)(M-2)}.
\]
Thus,
\[
(M-2)\int_{-\infty}^{t_0} f_{\bar
X}(z)\biggl(\int_{z}^{z+\delta}f_{\bar X}(t)\,\mathrm{d}\mu(t)\biggr)\,\mathrm
{d}\mu
(z)\leq
\frac{\delta\sqrt{\log(M-1)}}{M-1}.
\]
Note that if $t_0\geq\gamma_1$, then our claim follows. Indeed,
because $\Pro[\min_{2\leq j\leq M}\bar
X_j\leq\gamma_1] = 1-n^{-1}$, we have
\[
\Pro[\Omega_2] \geq\frac{1}{M-1}\biggl(1-\frac{1}{n}\biggr) - \delta\frac
{\sqrt
{\log(M-1)}}{M-1}.
\]
Otherwise, we split the interval $(-\infty,\gamma_1]=(-\infty
,t_0)\cup
[t_0,\gamma_1]$, and to
upper bound $T_2$, it remains to control the integral on the second
interval $[t_0,\gamma_1]$.

Recall that $W \in L_{\psi_1}$ and that $\E
W^3=0$. Therefore, by Proposition~\ref{prop:Petrov}, it is evident that
if $z$ and $\delta$ satisfy that $ z\leq z+\delta\leq0$ and
$|z|,|z+\delta|\leq B_0n^{1/6}$, then
%
\begin{eqnarray}\label{eq:Delta-part}
\int_{z}^{z+\delta}f_{\bar X}(t)\,\mathrm{d}\mu(t)&=&\Pro
[z\leq\bar X\leq z+\delta]
\nonumber
\\[-9pt]
\\[-9pt]
\nonumber
&\leq&\Pro[z\leq g\leq
z+\delta]+\frac{B_1}{\sqrt{n}}\exp(-z^2/2),
\end{eqnarray}
where $B_0$ and $B_1$ are constants that depend only on $\|W\|_{\psi
_1}$. In addition, for every $z\leq0$,
%
%
\begin{equation}\label{eq:Upper-Bound-Gaussian1}
\Pro[z\leq g\leq
z+\delta]\leq\frac{1}{\sqrt{2\uppi}}\exp(-z^2/2)
\int_0^\delta\exp(-zt)\,\mathrm{d}t\leq\frac{\delta}{\sqrt
{2\uppi}}\exp(-z^2/2).
\end{equation}

If $2M^2\exp(-B_0^2n^{1/3}/\|W\|_{\psi_2}^2)<\delta\leq1$, then
$|t_0|\leq B_0 n^{1/6}$. Combining \eqref{eq:Delta-part}
and \eqref{eq:Upper-Bound-Gaussian1} with the definition of $T_2$, we have
\begin{eqnarray*} 
&&(M-2)\int_{t_0}^{\gamma_1} f_{\bar X}(z)\biggl(\int
_{z}^{z+\delta
}f_{\bar X}(t)\,\mathrm{d}\mu(t)\biggr)\,\mathrm{d}\mu(z)
\\
&&\quad\leq(M-2)\biggl(\frac{B_1}{\sqrt{n}}+\frac{\delta}{\sqrt
{2\uppi
}}\biggr)\int_{t_0}^{\gamma_1}f_{\bar X}(z)
\exp(-z^2/2)\,\mathrm{d}\mu(z)
\\
&&\quad\leq(M-2)\biggl(\frac{B_1}{\sqrt{n}}+\frac{\delta}{\sqrt
{2\uppi}}\biggr)
\exp(-\gamma_1^2/2)\Pro[\bar X\leq\gamma_1]
\\
&&\quad\leq
(M-2)\biggl(\frac{B_1}{\sqrt{n}}+\frac{\delta}{\sqrt{2\uppi}}\biggr)
\exp(-\gamma_1^2/2)\frac{\log n}{M-1},
\end{eqnarray*}
where the last inequality follows from \eqref{eq:Def-gamma1-2}.
By Lemma~\ref{lemma:gamma-1}, and since $M \lesssim\sqrt{n}\log
n$,
\begin{eqnarray*}
&&(M-2)\int_{t_0}^{\gamma_1} f_{\bar
X}(z)\biggl(\int_{z}^{z+\delta}f_{\bar X}(t)\,\mathrm{d}\mu(t)\biggr)\,\mathrm
{d}\mu(z)
\\
&&\quad\leq c\biggl(\frac{1}{\sqrt{n}}+\delta\biggr)\biggl(\frac{\log
n}{M}\biggr)(\log n)\sqrt{\log M}
\end{eqnarray*}
for some constant $c=c(\beta)$, from which our claim follows.
\end{pf}

We next describe the construction needed for the proof of Theorem~\ref{thB}.
Let $(X,Y)$ and $F=\{f_1,\ldots,f_M\}$ be defined by
\begin{eqnarray*} \label{eq:counter-exemple}
Y & = &0,
\\
f_1(X) & =&(12)^{1/4}\cU_1,
\\
f_j(X)& =&(12)^{1/4}(\cU_j+\lambda) \qquad \mbox{for every }
2 \leq j \leq M,
\end{eqnarray*}
where $\cU_1,\ldots,\cU_M$ are $M$ independent random variables with
density $u\longmapsto
2(u+\lambda)\1_{[-\lambda,1-\lambda]}(u)$ for
$0<\lambda<1/2$ to be fixed later. Note that for this
choice of density function, $(\cU_1+\lambda)^2$ is uniformly
distributed on $[0,1]$, and the best element in $F$ with respect to the
quadratic risk is $f_1$.

Let $(\cU_j^{(i)}\dvt j=1,\ldots,M,i=1,\ldots,n)$ be a family of
independent random variables distributed as $\cU_1$. Thus, for
every $1 \leq i \leq n$, $f_j(X_i)=(12)^{1/4}(\cU_j^{(i)}+\lambda)$
for every $ 2\leq j \leq M$ and $f_1(X_i)=(12)^{1/4}\cU_1^{(i)}$.
For every $1 \leq j \leq M$, set
\[
\bar{R}_j=\sqrt{\frac{12}{n}}\Biggl(\sum_{i=1}^n\bigl(\cU_j^{(i)}+\lambda
\bigr)^2-\E\bigl(\cU
_j^{(i)}+\lambda\bigr)^2\Biggr),
\]
and observe that if
$W=\sqrt{12}(({\cU}+\lambda)^2-\E(\cU+\lambda)^2)$, then $W$
is a mean 0, variance $1$ random variable that is absolutely
continuous with respect to the Lebesgue measure and $W \in L_{\psi_2}$
and satisfies that $\E W^3=0$. These properties allow us to apply
Proposition~\ref{prop:Shahar-Proposition} to the random variables
$\bar
{R}_1,\ldots,\bar{R}_M$.

Let $0<\rho<1$ (to be named later), and set
\[
\label{eq:gamma-function}
\xi(\bar R_1)=\bar R_1+\frac{T}{\sqrt{n}}\log\biggl[
\frac{\rho}{2(1-\rho)}\biggr]-\sqrt{12}\lambda(2-\lambda)\sqrt{n},
\]
and
\[
\label{eq:delta}
\delta=\frac{-T}{\sqrt{n}}\log\biggl[\frac{\rho}{2(M-2)(1-\rho)}\biggr].
\]
Consider the system of inequalities
{\renewcommand{\theequation}{$C_j$}
\begin{equation}
\label{eq:System-Cj}
\cases{
\bar R_j\leq\xi(\bar R_1),\vspace*{2pt}\cr \bar R_k-\bar R_j\geq\delta\qquad \mbox{for
every } k\neq1,j,}
\end{equation}
}
\hspace*{-2pt}and recall that for each $j=1,\ldots,M$ $\hat{\theta}_j$ denotes the
weight of $f_j$ in the AEW procedure.

\begin{Proposition}\label{prop:Large-Weights}
There exist absolute constants $c_1$ and $c_2$ for which the following
holds. Let $0<\rho<1/2$ and $2 \leq j \leq M$. If the system (\ref{eq:System-Cj}) is
satisfied, then
\[
\hat\theta_j\geq1-\rho.
\]
Moreover, if $\rho\leq c_1 \lambda$, then the quadratic risk of the
function produced by the AEW procedure satisfies
\[
R(\tilde{f}^{\mathrm{AEW}})\geq\min_{f\in
F}R(f)+c_2\lambda.
\]
\end{Proposition}

\begin{pf}
Let $2 \leq j \leq M$, and assume that (\ref{eq:System-Cj}) is satisfied. Recall
that $R_n(f)$ is the empirical risk of $f$, and note that for any
$k\in\{2,\ldots,M\} \setminus\{j\}$,
\renewcommand{\theequation}{\arabic{section}.\arabic{equation}}
\setcounter{equation}{5}
\begin{eqnarray}
\label{eq:empirical-risk-kj}
\nonumber R_n(f_k)-R_n(f_j)&=&\frac{1}{n}\sum_{i=1}^n [f_k(X_i)^2-f_j(X_i)^2]
=\frac{\bar R_k-\bar R_j}{\sqrt{n}}
\nonumber
\\[-8pt]
\\[-8pt]
\nonumber
&\geq&\frac{\delta}{\sqrt{n}}
=\frac{-T}{n}\log\biggl[\frac{\rho}{2(M-2)(1-\rho)}\biggr].
\end{eqnarray}
In addition, since $\cU_1^{(i)}\leq1-\lambda$ almost surely for any $1
\leq i \leq
n$,
%
\begin{eqnarray} \label{eq:empirical-risk-1j}
R_n(f_1)-R_n(f_j)&=&
\frac{1}{n}\sum_{i=1}^n [f_1(X_i)^2-f_j(X_i)^2]
\nonumber\\
&=&\frac{\bar R_1-\bar R_j}{\sqrt{n}}-\sqrt{12}\Biggl(\lambda^2+\frac
{2\lambda
}{n}\sum_{i=1}^n\cU_1^{(i)}\Biggr)
\\
&\geq&\frac{\bar R_1-\xi(\bar
R_1)}{\sqrt{n}}-\sqrt{12}\lambda(2-\lambda)\geq
\frac{-T}{n}\log\biggl[\frac{\rho}{2(1-\rho)}\biggr].\nonumber
\end{eqnarray}
Combining \eqref{eq:empirical-risk-kj} and
\eqref{eq:empirical-risk-1j}, it is evident that
\begin{eqnarray*}
\hat\theta_j&=&\frac{1}{\sum_{k=1}^M
\exp[({-n}/{T})(R_n(f_k)-R_n(f_j))]}\\
&\geq&
\frac{1}{1+(M-2){\rho}/(2(M-2)(1-\rho))+{\rho}/(2(1-\rho))}=
1-\rho.
\end{eqnarray*}
Since the functions $f_1,\ldots,f_M$ are independent in $L_2(X)$ and
$\E f_j \geq0$,
\begin{eqnarray*}
R(\tilde{ f}^{\mathrm{AEW}})&=&\E\Biggl(\sum_{j=1}^M\hat\theta_j
f_j(X)\Biggr)^2\\
&=& (\hat{\theta}_j)^2 \E f_j^2 + \sum_{\ell\not= j}
(\hat{\theta}_\ell)^2 \E f_\ell^2 + \sum_{\ell\not= j}
\hat{\theta}_j \hat{\theta}_\ell\E f_j f_\ell
\geq(\hat{\theta}_j)^2 \E f_j^2,
\end{eqnarray*}
and there is an absolute constant $c_0$ for which $\E f_j^2 \geq
\E f_1^2 + c_0 \lambda$. Thus,
\[
(\hat{\theta}_j)^2 \E f_j^2 - \E f_1^2 \geq(1-\rho)(\E f_1^2 + c_0
\lambda) -
\E f_1^2 \geq c_2 \lambda,
\]
provided that $\rho\leq c_1 \lambda$, giving
\[
R(\tilde{ f}^{\mathrm{AEW}}) \geq\E f_1^2 + c_2 \lambda= \min_{f
\in
F} R(f)+c_2\lambda,
\]
as claimed.
\end{pf}

Next, we formulate a general statement, from which Theorem~\ref{thB} follows
immediately.

\begin{Theorem}
There exist absolute constants $c_i,i=0,\ldots,5$ and an integer $n_0$
for which
the following holds. For any $n \geq n_0$, $1 \leq\kappa\leq c_0\sqrt
{n \log
n}$, $0<T \leq1$, and $c_1T/\sqrt{n \log n}<\varepsilon<1/8$, let
$M=\lceil c_2\sqrt{n\log n}\rceil$, $\lambda=c_3 \varepsilon\sqrt
{(\log
n)/n}$, and
$\rho=n^{-\varepsilon\kappa/T}$. Set $F$ to be the class of
functions defined above with those parameters.
Then, with probability at least
\[
1-c_4(\varepsilon\kappa+T+1)\bigl((\log^3
n)/n\bigr)^{(1-2\varepsilon)^2/2},
\]
there exists $j\geq2$ such that
\[
\hat\theta_j\geq1-\frac{1}{n^{\varepsilon\kappa/T}}.
\]
In particular, with the same probability and if $0\leq
T<\min\{1,2\varepsilon\kappa\}$,
\[
R(\tilde{f}^{\mathrm{AEW}})\geq\min_{f\in F}R(f)+c_5\varepsilon
\sqrt
{\frac{\log
M}{n}}.
\]
\end{Theorem}

\begin{pf}
Set
\[
\Pro_0=\Pro[\exists j\in\{2,\ldots,M\} \mbox{ such that } \hat
\theta_j
\geq1-\rho],
\]
and, by Proposition~\ref{prop:Large-Weights},
\[
\Pro_0 \geq\Pro[\exists j \in\{2,\ldots,M\} \mbox{ for which (\ref{eq:System-Cj}) is satisfied}]=\Pro_1.
\]
Let $\gamma_1=\gamma_1(M-1)$ be defined by $\Pro[\min_{2\leq
j\leq M}\bar R_j\leq\gamma_1]=1-n^{-1}$, and observe that
$\gamma_1$ is well defined and satisfies all three parts of Lemma
\ref{lemma:gamma-1} for $\ell=M-1$. Set \mbox{$\Omega_0=\{\xi(\bar
R_1)\geq\gamma_1\}$},
\[
A=\bigl\{ \exists j \in\{2,\ldots,M\} \dvt\bar{R}_j \leq
\xi(\bar{R}_1), \mbox{ and } \bar{R}_k -\bar{R}_j \geq\delta{\mbox{
for every }} k \not= 1,j \bigr\}
\]
and
\[
B=\bigl\{\exists j \in\{2,\ldots,M\} \dvt\bar R_j\leq\gamma_1
\mbox{ and } \bar R_k-\bar R_j\geq\delta\mbox{{ for every }}
k\neq1,j \bigr\}.
\]
Since the functions $\bar R_j,j=1,\ldots,M$ are independent, we have
\[
\label{eq:P0}
\Pro_1 \geq\E_{\bar{R}_1} [ \Pro[A|\bar{R}_1]
\1_{\Omega_0} ] \geq\Pro[B] \Pro[\Omega_0].
\]
Applying Proposition~\ref{prop:Shahar-Proposition}, we then have
\[
\label{eq:Proba1}
\Pro[B] \geq1-\frac{1}{n} -c_2\biggl(\frac{1}{\sqrt{n}}+\delta\biggr)(\log
n)^2\sqrt{\log M},
\]
provided that $c_3\log n\leq M\leq c_4\sqrt{n}(\log n)$.

To lower bound $\Pro[\Omega_0]$, note that
\[
\label{eq:Proba-Omega0}
\Pro[\Omega_0]=\Pro\biggl[\bar R_1\geq
\gamma_1-\frac{T}{\sqrt{n}}\log\biggl(\frac{\rho}{2(1-\rho)}\biggr)
+\sqrt{12}\lambda(2-\lambda)\sqrt{n}\biggr].
\]
Fix $0<\varepsilon<1/8$ and assume that $\lambda,\rho$ and $T$ are
such that
%
%
\begin{equation} \label{eq:def-theta-rho}
\sqrt{12}\lambda(2-\lambda)\sqrt{n}\leq-\varepsilon\gamma_1
\quad\mbox{and}\quad
{-}\frac{T}{\sqrt{n}}\log\biggl(\frac{\rho}{2(1-\rho)}\biggr)\leq-\varepsilon
\gamma_1.
\end{equation}
By the Berry--Ess{\'e}en theorem and
\eqref{eq:Gaussian-tail-estimate},
%
\begin{eqnarray*}\label{eq:Proba-Omega0-2}
\Pro[\Omega_0]&\geq&\Pro[\bar
R_1\geq(1-2\varepsilon)\gamma_1] =1-\Pro[\bar
R_1<(1-2\varepsilon)\gamma_1] \\
&\geq& 1-\Pro[g\leq
(1-2\varepsilon)\gamma_1]-\frac{2\beta(W)}{\sqrt{n}}
\\
& \geq&
1-\frac{1}{\sqrt{2\uppi}(1-2\varepsilon)|\gamma_1|}\exp
\bigl(-(1-2\varepsilon
)^2\gamma_1^2/2\bigr)
-\frac{2A}{\sqrt{n}},
\end{eqnarray*}
and by Lemma~\ref{lemma:gamma-1},
\[
\exp\bigl(-(1-2\varepsilon)^2\gamma_1^2/2\bigr) \leq c_5\biggl(\frac{\log
n}{M-1}{\log
^{1/2}\biggl(\frac{c_5M}{\log n}\biggr)}\biggr)^{(1-2\eps)^2}.
\]
Therefore,
\[
\Pro_0\geq\biggl(1-\frac{1}{n}-c_2\biggl(\frac{1}{\sqrt{n}}+\delta\biggr)
(\log n)^2\sqrt{\log M}\biggr)
\cdot\biggl(1-c_5\biggl(\frac{\log^3
n}{M}\biggr)^{(1-2\varepsilon)^2}\biggr),
\]
provided that $c_2\log n\leq M\leq c_3\sqrt{n\log n}$.

To complete the proof, we need to chose $\lambda$ and $\rho$
for which \eqref{eq:def-theta-rho} holds. By Lemma
\ref{lemma:gamma-1},
\[
|\gamma_1| \gtrsim\log^{1/2}\biggl(\frac{M}{\log n}\biggr),
\]
and thus \eqref{eq:def-theta-rho} holds for $\lambda$ and $\rho$ for
which
\[
\lambda\leq c_8\varepsilon\biggl[\frac{1}{n}\log\biggl(\frac{M}{\log
n}\biggr)\biggr]^{1/2}\quad \mbox{and}\quad \rho\geq
2\exp\biggl[\frac{-c_9\varepsilon\sqrt{n}}{T}\log^{1/2}\biggl(\frac{M}{\log
n}\biggr)\biggr].
\]
In particular, when we take $ M \sim\sqrt{n\log n}$, $\lambda\sim
\varepsilon((\log M)/n)^{1/2}$, and $\rho=n^{-\varepsilon
\kappa/T}$, $\rho$ satisfies the required condition as long as
$\varepsilon\gtrsim T/\sqrt{n \log n}$ and $\kappa\lesssim\sqrt
{n/\log
n}$, as assumed. Moreover,
\[
\delta\lesssim(\varepsilon\kappa+ T)\frac{\log{n}}{\sqrt{n}},
\]
implying that
\[
\Pro_0\geq1-c_8(\varepsilon\kappa+ T +1)\biggl(\frac{\log^3
n}{n}\biggr)^{{(1-2\varepsilon)^2}/{2}}.
\]
The lower bound on the risk of the AEW procedure now follows from
Proposition~\ref{prop:Large-Weights}.
\end{pf}

\section{\texorpdfstring{Proof of Theorem \protect\ref{thC}}{Proof of Theorem C}} \label{sec:proof-of-C}
In this section we prove Theorem~\ref{thC}, which we reformulate
below. From here on, we assume that the dictionary $F$ is
finite, consisting of $M$ functions, and that the functions are
indexed according to their risk in an increasing order. Thus,
$f_1=f_F^*$. In addition, we denote $\mathcal{ L}_f(\cdot) =
Q(\cdot,f)-Q(\cdot,f_1)$, and thus $R(f)-R(f_1)=\E\mathcal{ L}_f$.

For every $r>0$, recall that
\begin{eqnarray*}
\psi(r) &= & \log(|\{ f \in F \dvt\E\mathcal{ L}_f \leq r \}|+1)
\\
&&{}+ \sum_{j=1}^\infty2^{-j} \log(|\{ f \in F \dvt 2^{j-1} r < \E
\mathcal{ L}_f \leq2^j r \}|+1),
\end{eqnarray*}
which serves as a measure of complexity for the class $F$.

The first component needed in the proof of Theorem~\ref{thC} is the level
$\lambda(x)$ with the following property: with
probability at least $1-2\exp(-x)$, $R_n(f_j)-R_n(f_1)$ is
equivalent to $R(f_j)-R(f_1)$ if $R(f_j)-R(f_1) \geq\lambda(x)$.
This ``isomorphism'' constant was introduced by~\cite{MR2240689}. To
formulate the exact properties that we need, first recall the following
definitions and notation.

If $G=\mathcal{ L}_F$ is the excess loss functions class $\{\cL
_f\dvt
f\in
F\}$, then let $\operatorname{star}(G,0)=\{\theta g \dvt 0 \leq
\theta\leq1, g
\in G\}$ is the star-shaped hull of $G$ and $0$. Set $G_r=\operatorname
{star}(G,0) \cap\{g \dvt \E g =r\}$, that is, the set of functions in
the star-shaped hull of $\mathcal{ L}_F$ and $0$, with expectation
$r$. Let
\[
r^*=\inf\Bigl\{r\dvt\E\sup_{g \in G_r} |P_n g - Pg| \leq r/2\Bigr\},
\]
where, as always, $P_n$ denotes the empirical mean and $P$ is the mean
according to the underlying probability measure of $Z$.

\begin{Theorem}[(\cite{MR2240689})] \label{thm:isomorphic-coordinate-proj}
 There exists an absolute constant $c$ for which the
following holds. Let $F$ be a class of functions bounded by $b$,
such that $\cL_F$ is a $(1,B)$-Bernstein class. For every $x>0$ and
an integer $n$, let
%
%
\begin{equation}
\label{eq:lambda}
\lambda(x)=c \max\biggl\{r^*, (b+B)\frac{x}{n} \biggr\}.
\end{equation}
Then, with probability at least $1-2\exp(-x)$, for every $f \in F$
with $R(f)-R(f^*_F) \geq\lambda(x)$,
\[
R_n(f)-R_n(f_F^*) \geq\frac{1}{2}\bigl(R(f)-R(f_F^*)\bigr).
\]
\end{Theorem}

Let $\rho=\kappa_1 (B +b)/n$, where $\kappa_1$ is an absolute
constant to be named later. Recall that functions in $F$ are indexed
according to their risk in an increasing order. Let $J_-(x)=\{j \dvt
R(f_j)-R(f_1) \leq\lambda(x)\}$, and set $J_+(x)$ as its
complement. Define the sets $J_{+,0}=\{j \in J_+(x)\dvt R(f_j)-R(f_1)
\leq\rho\}$ and, for $k \geq1$,
\[
J_{+,k} = \{j \in J_+(x)\dvt 2^{k-1} \rho< R(f_j)-R(f_1) \leq2^{k}
\rho\}.
\]
(Note that some of the sets $J_{+,k}$ may be empty.) Set
\[
k_0 = \sup\{ k \geq0 \dvt 2^k \leq\log(|J_{+,k}|+1) \},
\]
and let $I=J_{-} \cup\bigcup_{k \leq k_0} J_{+,k}$.

From Theorem~\ref{thm:isomorphic-coordinate-proj}, it follows that
for every $k \geq0$ and every $j \in J_{+,k}$, $R_n(f_j)-R_n(f_F^*)
\geq\frac{1}{2}(R(f_j)-R(f_F^*))$. This is because
$R(f_j)-R(f_F^*)\geq\lambda(x)$ by the definition of $J_+(x)$, and
$J_+(x)\supset J_{+,k}$.

The key factor in the proof of Theorem~\ref{thC} is Theorem
\ref{thm:optimal-agg}.

\begin{Theorem} \label{thm:optimal-agg}
There exist absolute constants $c_1$ and $c_2$ for which the
following holds. Let $F$ be a class of functions bounded by $b$,
such that $\cL_F$ is a $(1,B)$-Bernstein class with respect to a
convex risk function $R$. Then, with probability at least
$1-2\exp(-x)$, if $\tilde{f}^{\mathrm{AEW}}$ is produced by the AEW
algorithm and
$T \leq c_1(b+B)$, then
%
%
\begin{equation} \label{eq:key-estimate}
R(\tilde{f}^{\mathrm{AEW}})-R(f_F^*) \leq c_2 \biggl(\lambda(x)+(b+B)
\frac{2^{k_0}}{n}\biggr),
\end{equation}
where $\lambda(x)$ is as defined in~(\ref{eq:lambda}).
\end{Theorem}

\begin{pf} Let $(\hat{\theta}_j)_{j=1}^M$ be the weights of the AEW
algorithm, and set $\tilde{f}^{\mathrm{AEW}} = \sum_{j=1}^M \hat
{\theta
}_j f_j$ to
be the aggregate function. Because $R$ is a convex function,
\[
R\Biggl(\sum_{j=1}^M \hat{\theta}_j f_j \Biggr) -R(f_1) \leq
\sum_{j=1}^M \hat{\theta}_j \bigl(R(f_j)-R(f_1)\bigr).
\]

Note that for every $j \in I$, $R(f_j)-R(f_1) \leq\lambda(x)
+2^{k_0} \rho= \lambda(x) + \kappa_1 2^{k_0}(b+B)/n$. In
particular, because $\sum_{j=1}^M \hat{\theta}_j = 1$,
\[
\sum_{j \in I} \hat{\theta}_j \bigl(R(f_j)-R(f_1)\bigr) \leq\lambda(x) +
\kappa_1 2^{k_0}(b+B)/n.
\]
On the other hand, with probability at least $1-2\exp(-x)$, for
every $k>k_0$ and every $j \in J_{+,k}$,
\[
R_n(f_j) - R_n(f_1) \geq\bigl(R(f_j)-R(f_1)\bigr)/2.
\]
Applying the definition of the weights in the AEW algorithm and
given that $\hat{\theta}_1 \leq1$,
\begin{eqnarray*}
\sum_{j \in I^c} \hat{\theta}_j \bigl(R(f_j)-R(f_1)\bigr)& =&
\hat{\theta}_1 \sum_{j \in I^c}
\frac{\hat{\theta}_j}{\hat{\theta}_1} \bigl(R(f_j)-R(f_1)\bigr)
\\
&\leq& \sum_{j \in I^c} \exp\biggl(-\frac{n}{T}
\bigl(R_n(f_j)-R_n(f_1)\bigr) \biggr) \bigl(R(f_j)-R(f_1)\bigr)
\\
&\leq& \sum_{k > k_0 } \sum_{j \in J_{+,k}} \exp\biggl(-\frac{n}{2T}
\bigl({R(f_j)-R(f_1)}\bigr)\biggr) \bigl(R(f_j)-R(f_1)\bigr) =
(\star).
\end{eqnarray*}
From the definition of $k_0$, it is evident that for every $k > k_0$,
$2^k \geq\log|J_{+,k}|$, and thus if $T \leq c_1 \max\{b,B\}$ and
$\kappa_1$ is sufficiently large, then
\[
(\star)\leq\sum_{k > k_0 } \exp\biggl(\log|J_{+,k}| - \frac{n}{2T}
2^{k-1} \rho\biggr)2^{k}\rho\leq\sum_{k>k_0} \exp\biggl(-c_2
\frac{n}{T}2^k\rho\biggr) 2^{k} \rho\leq c_3 \frac{T}{n}.
\]
Indeed, this follows because for that choice of $T$,
$(n/T)2^{k_0}\rho\geq c_4$, with $c_4$ an absolute constant.

Thus, with probability at least $1-2\exp(-x)$,
\[
R(\tilde{f}) -R(f_1) \leq\lambda(x) + \kappa_1 2^{k_0}(b+B)/n + c_3
\frac{T}{n} \leq\lambda(x) + c_5 2^{k_0} \frac{b+B}{n},
\]
as claimed.
\end{pf}

The next step in the proof of Theorem~\ref{thC} requires several simple
facts regarding the empirical process indexed by a localization of
the star-shaped hull of a Bernstein class. First, it is simple to
verify that the star-shaped hull of a
$(1,B)$-Bernstein class is a $(1,B)$-Bernstein class as well.
Second, if $G =\operatorname{star}(\mathcal{ L}_F,0)$ and $G_r=\{h
\in
G \dvt\E h
=r\}$, then
\[
G_r = \bigcup_{j \geq1} \biggl\{ \frac{r \mathcal{ L}_f }{\E\mathcal{ L}_f}
\dvt f\in F, 2^{j-1} r \leq\E\mathcal{ L}_f \leq2^j r \biggr\} \equiv
\bigcup_{j \geq1} H_{r,j}.
\]
In particular,
\[
\E\sup_{h \in G_r} \Biggl| \frac{1}{n} \sum_{i=1}^n h(Z_i) - \E h
\Biggr| \leq\sum_{i=1}^\infty\E\sup_{h \in H_{r,j}} \Biggl|
\frac{1}{n} \sum_{i=1}^n h(Z_i) - \E h \Biggr|.
\]

\begin{Lemma} \label{lemma:localized-estimates}
There exists an absolute constant $c$ for which the following holds.
If $\cL_F$ is a $(1,B)$-Bernstein class with respect to $Z$, then for every
$r$ and $j \geq1$,
\[
\E\sup_{h \in H_{r,j}} |P_n h - P h | \leq c
\max\Biggl\{\frac{b2^{-j} \log(|H_{r,j}|+1)}{n}, \sqrt{\frac{\log
(|H_{r,j}|+1)}{n}} \sqrt{rB2^{-j}}\Biggr\}.
\]
\end{Lemma}

\begin{pf} Fix $r>0$ and $j \geq1$, and let
\[
D=\sup_{h \in H_{r,j}} \Biggl(\frac{1}{n}\sum_{i=1}^n
h^2(Z_i)\Biggr)^{1/2}.
\]
Note that every $h \in H_{r,j}$ satisfies that $h=r \mathcal{ L}_f/\E
\mathcal{ L}_f$ for some $f\in F$, and for which $\E\mathcal{ L}_f
\geq
r2^{j-1}$. Therefore, using the Bernstein condition on $\cL_F$,
\[
\E h^2 = r^2 \frac{\E(\mathcal{ L}_f)^2}{(\E\mathcal{ L}_f)^2}
\leq
rB2^{-j+1}.
\]
Moreover, $\|h\|_\infty\leq(r/ \E\mathcal{ L}_f) \|\mathcal{
L}_f\|_\infty\leq b2^{-j+1}$. Thus, by the Gin\'{e}--Zinn
symmetrization theorem and a contraction argument
(see, e.g.,~\cite{MR757767} and~\cite{LT:91}),
\begin{eqnarray*}
\E D^2 &\leq& \E\sup_{h \in H_{r,j}}\Biggl |\frac{1}{n}\sum_{i=1}^n
h^2(Z_i) - \E h^2 \Biggr| + rB2^{-j+1}
\\
&\leq& \frac{2}{\sqrt{n}}\E_Z \E_\eps\sup_{h \in H_{r,j}}
\Biggl|\frac{1}{\sqrt{n}} \sum_{i=1}^n \eps_i h^2(Z_i) \Biggr| +
rB2^{-j+1}
\\
&\leq& \frac{b2^{-j+2}}{\sqrt{n}} \E_Z \E_\eps\sup_{h \in H_{r,j}}
\Biggl|\frac{1}{\sqrt{n}} \sum_{i=1}^n \eps_i h(Z_i)\Biggr | +
rB2^{-j+1}
\\
&\leq& \frac{c_0rb2^{-j+2}}{\sqrt{n}} \sqrt{\log(|H_{r,j}|+1)} \E D
+ rB2^{-j+1},
\end{eqnarray*}
where the last inequality is evident by the sub-Gaussian properties
of the Rademacher process (cf.~\cite{LT:91}). Since $\E D \leq(\E
D^2)^{1/2}$, it follows that
\[
\E D^2 \leq c_0 b2^{-j+2}\sqrt{\frac{\log(|H_{r,j}|+1)}{n}} (\E
D^2)^{1/2} + rB2^{-j+1},
\]
implying that
\[
\E D^2 \leq c_1\max\biggl\{b^22^{-2j} \frac{\log(|H_{r,j}|+1)}{n},
rB2^{-j}\biggr\}.
\]
Thus, again using a symmetrization argument and the sub-Gaussian
properties of the Rademacher process, we have
\begin{eqnarray*}
\hspace*{-20pt} \E\sup_{h \in H_{r,j}} \Biggl|\frac{1}{n} \sum_{i=1}^n h(Z_i) - \E
h \Biggr| &\leq&\frac{c_2}{\sqrt{n}} \sqrt{\log(|H_{r,j}|+1)} \E D
\\
&\leq& c_3 \max\Biggl\{\frac{b2^{-j} \log(|H_{r,j}|+1)}{n},
\sqrt{\frac{\log(|H_{r,j}|+1)}{n}} \sqrt{rB2^{-j}}\Biggr\}.
\end{eqnarray*}
\upqed\end{pf}

\begin{Corollary} \label{cor:est-for-r^*}
There exist absolute constants $c_1$ and $c_2$ for which the
following holds. Let $F$ be a finite class consisting of $M$
functions bounded by $b$, such that the excess loss class $\cL_F$ is
a $(1,B)$-Bernstein class. If we set $\theta=c_1(b+B)(\log M)/n$,
then
\[
r^* \leq c_2\biggl(\frac{b+B}{n}\biggr)\psi(\theta).
\]
\end{Corollary}

\begin{pf} Observe that for every $r>0$,
\begin{eqnarray*}
 &&\E\sup_{h \in G_r}\Biggl |\frac{1}{n}\sum_{i=1}^n h(Z_i) - \E h
\Biggr| \\
&&\quad\leq\sum_{j \geq1} \E\sup_{h \in H_{r,j}}
\Biggl|\frac{1}{n} \sum_{i=1}^n h(Z_i) - \E h \Biggr|
\\
&&\quad\leq c_1 \max\Biggl\{ \frac{b}{n} \sum_{j \geq1} 2^{-j} \log
(|H_{r,j}|+1), \sqrt{\frac{Br}{n}} \sum_{j \geq1} 2^{-j/2}
\sqrt{\log(|H_{r,j}|+1)} \Biggr\}
\\
&&\quad\leq c_1 \frac{b}{n} \biggl(\log(|H_{r,0}|+1)+ \sum_{j \geq1}
2^{-j} \log(|H_{r,j}|+1) \biggr)
\\
&&\qquad{}+  c_1\sqrt{\frac{Br}{n}} \biggl(\sqrt{\log(|H_{r,0}|+1)}+\sum_{j
\geq1} 2^{-j/2} \sqrt{\log(|H_{r,j}|+1)}\biggr)
\\
&&\quad\equiv u(r),
\end{eqnarray*}
where we define $H_{r,0}=\{(r\cL_f)/(\E\cL_f)\dvt f\in F, \E
\cL_f\leq r\}$. Let $\bar{r}=\inf\{r \dvt u(r) \leq r/2\}$. Since
$|H_{r,j}| \leq M$ for every $j \geq0$, we have
\[
u(r) \leq c_2 \max\Biggl\{b\frac{\log M}{n}, \sqrt{\frac{rB \log
M}{n}} \Biggr\},
\]
and thus
\[
\bar{r} \leq c_3 (b+B)(\log M)/n = \theta.
\]
Moreover, the functions of $r$,
\[
\log(|H_{r,0}|+1)+\sum_{j \geq1} 2^{-j} \log(|H_{r,j}|+1),
\]
and
\[
\sqrt{\log(|H_{r,0}|+1)}+\sum_{j \geq1} 2^{-j/2} \sqrt{\log
(|H_{r,j}|+1)},
\]
are increasing, and thus for any $r \leq\theta$,
\begin{eqnarray*}
&& \frac{b}{n} \biggl(\log(|H_{r,0}|+1)+ \sum_{j \geq1} 2^{-j} \log
(|H_{r,j}|+1)\biggr)
\\
&&\quad\leq \frac{b}{n} \biggl(\log(|H_{\theta,0}|+1)+ \sum_{j \geq1}
2^{-j} \log(|H_{\theta,j}|+1)\biggr)
\end{eqnarray*}
and
\begin{eqnarray*}
&& \sqrt{\frac{Br}{n}} \biggl(\sqrt{\log(|H_{r,0}|+1)}+\sum_{j \geq
1} 2^{-j/2} \sqrt{\log(|H_{r,j}|+1)}\biggr)
\\
&&\quad\leq \sqrt{\frac{Br}{n}} \biggl(\sqrt{\log
(|H_{\theta,0}|+1)}+\sum_{j \geq1} 2^{-j/2} \sqrt{\log
(|H_{\theta,j}|+1)}\biggr).
\end{eqnarray*}
Thus, if we consider
\begin{eqnarray*}
r& = & c_3\frac{b}{n} \biggl(\log(|H_{\theta,0}|+1)+ \sum_{j \geq1}
2^{-j} \log(|H_{\theta,j}|+1) \biggr)
\\
&&{}+  c_3 \frac{B}{n}\biggl(\sqrt{\log(|H_{\theta,0}|+1)}+\sum_{j \geq
1} 2^{-j/2} \sqrt{\log(|H_{\theta,j}|+1)}\biggr)^2
\\
&\leq& c_4\biggl(\frac{b+B}{n}\biggr)\psi(\theta)
\end{eqnarray*}
for appropriate constants $c_3$ and $c_4$, then $r\leq\theta$.
Thus, $u(r) \leq r/2$ and, therefore,
\[
\bar{r} \leq
c_4\biggl(\frac{b+B}{n}\biggr)\psi(\theta).
\]
Finally, because
\[
\E
\sup_{h \in G_r} |P_nh-Ph| \leq u(r)
\]
and $r^*=\inf\{r \dvt\E
\sup_{g \in G_r} |P_n g - Pg| \leq r/2\}$, we have $r^* \leq\bar{r}$.
\end{pf}

\begin{pf*}{Proof of Theorem \protect\ref{thC}} The proof of Theorem~\ref{thC} follows
from estimates of $\lambda(x)$ and $2^{k_0}$. From Corollary \ref
{cor:est-for-r^*}, it is evident that
\[
\lambda(x) \leq c_1 \max\biggl\{ \biggl(\frac{b+B}{n}\biggr)\psi\biggl(c_1(b+B)\frac{\log
M}{n}\biggr), (b+B)\frac{x}{n} \biggr\},
\]
where $c_1$ is an absolute constant to be identified later. (Note that
$\psi$ is an increasing function.)

Next, by the definition of $k_0$, $2^{k_0} \leq\log M$.
Therefore, using the notation of Theorem~\ref{thm:optimal-agg},
\[
\bigcup_{k \leq k_0} \{ f_j \dvt j \in J_{+,k} \} \subset
\biggl\{f_j \dvt R(f_j)-R(f_1) \leq\kappa_1(b+B)\frac{\log M}{n}
\biggr\}
\]
and, in particular,
%
\begin{eqnarray*}
2^{k_0} & \leq&\log\biggl(\biggl|\bigcup_{k \leq k_0} \{ f_j \dvtx j \in
J_{+,k} \} \biggr|+1\biggr)
\\
& \leq&\log\biggl(\biggl| \biggl\{f_j \dvt R(f_j)-R(f_1) \leq
\kappa_1(b+B)\frac{\log M}{n} \biggr\} \biggr|+1\biggr) \leq\log
(|H_{\theta,0}|+1),
\end{eqnarray*}
for an appropriate choice of constant $c_1$.

The second part of Theorem~\ref{thC} follows from a standard integration
argument.
\end{pf*}

\section*{Acknowledgements}
This article was written while G. Lecu\'{e} was visiting the
Department of Mathematics, Technion, and the Centre for Mathematics
and Its Applications, Australian National University. The
authors thank both of these institutions for their
hospitality. They also thank Pierre Alquier and Olivier Catoni for
useful discussions.
G. Lecu{}\'e was supported by French Agence Nationale de la
Recherche ANR Grant \textsc{``Prognostic''} ANR-09-JCJC-0101-01.
S. Mendelson was supported in part by the Centre for Mathematics and
its Applications, The Australian National University, Canberra, ACT
0200, Australia, by an Australian Research Council Discovery Grant
DP0559465, DP0986563 and by the European Community's Seventh
Framework Programme (FP7/2007-2013), ERC grant agreement 203134.

%

\printhistory

\end{document}